\documentclass[11pt,reqno]{amsart}

\numberwithin{equation}{section}
\newtheorem{thm}{Theorem}[section]

\newtheorem{cor}[thm]{Corollary}

\newtheorem{theorem}{Theorem}[section]
\newtheorem{lemma}[theorem]{Lemma}

\theoremstyle{definition}
\newtheorem{example}[theorem]{Example}

\newtheorem{remark}[theorem]{Remark}

\newtheorem*{acks}{Acknowledgements}

\newenvironment{romenumerate}{\begin{enumerate}% gives (i), (ii) etc.
 }{\end{enumerate}}

\newenvironment{thmxenumerate}{\begin{enumerate}% gives (i), (ii) etc.
\setlength{\leftmargin}{0pt}
\setlength{\itemindent}{0pt}
 }{\end{enumerate}}

\newcounter{thmenumerate}

\newcommand\pfitem[1]{\par(#1):}

\newcommand{\AAA}{A${}'$}

\newcommand{\condition}[1]{Condition~\textup{#1}}
\newcommand{\conditionbd}{\condition{B($\gd$)}}
\newtheorem*{conditionA}{Condition~A}
\newtheorem*{conditionAAA}{Condition~\AAA}
\newtheorem*{conditionBd}{Condition~B($\gd$)}

\newcommand{\convention}[1]{Convention \textup{#1}}
\newtheorem*{conventionC}{\convention C}

\newcommand{\refT}[1]{Theorem~\ref{#1}}
\newcommand{\refC}[1]{Corollary~\ref{#1}}
\newcommand{\refL}[1]{Lemma~\ref{#1}}
\newcommand{\refR}[1]{Remark~\ref{#1}}
\newcommand{\refS}[1]{Section~\ref{#1}}

\newcommand{\refE}[1]{Example~\ref{#1}}
\newcommand{\refand}[2]{\ref{#1} and~\ref{#2}}

\newcommand\marginal[1]{\marginpar{\raggedright\parindent=0pt\tiny #1}}
\newcommand{\ignore}[1]{}
\newcommand\REM[1]{\texttt{[#1]}\marginal{XXX}}

\newcommand\set[1]{\ensuremath{\{#1\}}}
\newcommand\bigset[1]{\ensuremath{\bigl\{#1\bigr\}}}

\newcommand\xpar[1]{(#1)}
\newcommand\bigpar[1]{\bigl(#1\bigr)}
\newcommand\Bigpar[1]{\Bigl(#1\Bigr)}

\newcommand\bbigpar[1]{\left(#1\right)}

\newcommand\Bigparq[1]{\Bigl[#1\Bigr]}

\newcommand\lrabs[1]{\left|#1\right|}
\def\rompar(#1){\textup(#1\textup)}    % usage: \rompar(...)
\newcommand\xfrac[2]{#1/#2}
\newcommand\parfrac[2]{\Bigpar{\frac{#1}{#2}}}
\newcommand\norm[1]{\ensuremath{\|#1\|}}

\newcommand\Norm[1]{\ensuremath{\left\|#1\right\|}}

\newcommand\ceil[1]{\lceil#1\rceil}

\newcommand\ttoo{\ensuremath{{t\to\infty}}}
\newcommand\xtoo{\ensuremath{{x\to\infty}}}

\def\[#1]{[\![#1]\!]}
\newcommand\ett{\mathbf 1}
\newcommand\etta[1]{\ett_{#1}}
\newcommand\ettax[1]{\ett_{\set{#1}}}
\newcommand\ind{\ettax}
\newcommand\Bigind[1]{\ett_{\bigset{#1}}}

\newcommand\iid{i.i.d.\spacefactor=1000}     %?????
\newcommand\ie{i.e.\spacefactor=1000}
\newcommand\eg{e.g.\spacefactor=1000}

\newcommand{\as}{a.s.\spacefactor=1000}

\newcommand\dto{\overset{\mathrm{d}}{\to}}

\newcommand\eqd{\overset{\mathrm{d}}{=}}
\renewcommand\={:=}

\newcommand\bbR{\mathbb R}
\newcommand\bbC{\mathbb C}

\newcommand\ii{\mathrm i}

\newcommand\E{\operatorname{\mathbb E{}}}
\renewcommand\P{\operatorname{\mathbb P{}}}
\newcommand\Var{\operatorname{Var}}
\newcommand\Cov{\operatorname{Cov}}

\newcommand\unif{\operatorname{U{}}}
\newcommand\Res{\operatorname{Res}}
\renewcommand\Re{\operatorname{Re}}
\renewcommand\Im{\operatorname{Im}}
\newcommand\dd{\,d}

\newcommand\intoo{\int_0^\infty}
\newcommand\intooo{\int_{-\infty}^\infty}
\newcommand\intoi{\int_0^1}
\newcommand\sumjb{\sum_{j=1}^b}
\newcommand\sumkb{\sum_{k=1}^b}
\newcommand\sumn{\sum_{n=0}^\infty}

\newcommand\ga{\alpha}
\newcommand\gb{\beta}
\newcommand\gd{\delta}
\newcommand\gD{\Delta}

\newcommand\gl{\lambda}
\newcommand\gL{\Lambda}
\newcommand\go{\omega}

\newcommand\gs{\sigma}
\newcommand\gss{\sigma^2}

\newcommand\eps{\varepsilon}

\newcommand\cB{\mathcal B}

\newcommand\cF{\mathcal F}

\newcommand\cH{\mathcal H}
\newcommand\cL{{\mathcal L}}
\newcommand\cM{\mathcal M}
\newcommand\cN{\mathcal N}

\newcommand\cT{{\mathcal T}}

\newcommand\hatf{\hat f}

\newcommand\bmin{\wedge}
\newcommand\bmax{\vee}

\newcommand\qi{^{-1}}
\newcommand\qh{^{1/2}}
\newcommand\qhi{^{-1/2}}

\newcommand{\nc}{\newcommand}
\nc{\Prob}{{\mathbb{P}}}
\nc{\Nset}{{\mathbb{N}}}
\nc{\Rset}{{\mathbb{R}}}
\nc{\Zset}{{\mathbb{Z}}}
\nc{\bes}{\begin{eqnarray*}}
\nc{\ees}{\end{eqnarray*}}
\nc{\be}{\begin{eqnarray}}
\nc{\ee}{\end{eqnarray}}

\nc{\op}{{\mathrm{op}}}
\nc{\Id}{{\mathrm{Id}}}
\nc{\deq}{{\;\stackrel{d}{=}}\;}
\nc{\vv}{\ensuremath{\mathbf{V}}}
\nc{\vvb}{\ensuremath{(V_1,\dots,V_b)}}
\nc{\tx}{\ensuremath{\cT(x)}}
\nc{\too}{\ensuremath{\cT(\infty)}}
\nc{\nx}{\ensuremath{N(x)}}
\nc{\nex}{\ensuremath{N_e(x)}}
\nc{\ny}{\ensuremath{N(y)}}
\nc\cBx{\cB^\ast}
\nc\jj[1]{j_1\cdots j_{#1}}

\nc\nxx{N_*}
\nc\nxt{\nxx(t)}
\nc\mxx{m_*}
\nc\mxt{\mxx(t)}
\nc\gssx{\gss_*}
\nc\gssxt{\gssx(t)}
\nc\tf{\tilde f}
\nc\tF{\tilde F}
\nc\tH{\widetilde H}
\nc\tmu{\widetilde \mu}
\nc\tnu{\widetilde \nu}
\nc\tmxx{\widetilde \mxx}
\nc\gai{\ga\qi}
\nc\epsi{\eps\qi}
\nc\epsix{\frac1{\eps}}
\nc\ff{f^*}
\nc\ooo{\ensuremath{[0,\infty)}}
\nc\he{H_\eps}
\nc\tHe{\widetilde\he}
\nc\mex{m_{*\eps}}
\nc\tmex{\widetilde\mex}
\nc\gsx{\gs_2}
\nc\hiu{1/2+\ii u}
\nc\hiuu{1/2-\ii u}
\nc\gDx{\Delta^{\!\ast}}
\nc\gDxoo{\gDx(\infty)}
\nc\taud{\tau(\gd)}

\nc\Nx{M}

\newcommand{\che}{characteristic equation}

\newcommand{\Renyi}{R\'enyi}

\nc\tgs{\tilde\gs}
\nc\gst{\gs_t}
\nc\gsst{\gss_t}
\nc\trt{T_r^{(t)}}
\nc\art{A_r^{(t)}}
\nc\bt{b^{(t)}}
\newcommand\lsto{\overset{\ell_s}{\longrightarrow}}
\newcommand\ax{A^\ast}
\newcommand\axr{\ax_r}
\nc\taue{\tau_\eps}
\newcommand\chd{\cH_\gd}
\newcommand\bchd{\overline{\cH_\gd}}
\newcommand\bz{\bar{z}}
\newcommand\taux{t}
\begin{document}

\title%[]
{The size of random fragmentation trees}

%\makeatletter
%\let\@fnsymbol\@arabic % \@roman, \@Roman, \@alph, \@Alph
%\makeatother
\nc\tss{\textsuperscript}

\newcommand\urladdrx[1]{\urladdr{\def~{$\sim$}#1}}

\author{Svante Janson}
\address{Department of Mathematics, Uppsala University, PO Box 480,
SE-751~06 Uppsala, Sweden}
\email{svante.janson@math.uu.se}
\urladdrx{http://www.math.uu.se/~svante/}

\author{Ralph Neininger\tss1}
\thanks{\tss1Research supported by an Emmy Noether fellowship of the DFG}
\address{Department of Mathematics and Computer Science,
J.W.~Goethe University, 60054 Frankfurt a.M.,
Germany} \email{neiningr@math.uni-frankfurt.de}
\urladdrx{http://www.math.uni-frankfurt.de/~neiningr/}

\date{September 13, 2006}

%\keywords{<keywords>}
%\subjclass{Primary: <subject>; Secondary: <subject>}

\begin{abstract}
We study a random fragmentation process and its associated random
tree.
The process has earlier been studied by Dean and Majumdar
\cite{DeanMaj},
who found a phase transition: the number of fragmentations is asymptotically
normal in some cases but not in others, depending on the position of
roots of a certain characteristic equation.
This parallels the behaviour of discrete analogues with various random
trees that have been studied in computer science.
We give rigorous proofs of this phase transition, and add further details.

The proof uses the contraction method. We extend some previous results for
recursive sequences of random variables to
families of random variables with a continuous parameter;
we believe that this extension has independent interest.
\end{abstract}

\maketitle
\renewcommand\tss[1]{}

\section{The problem and result}

Consider the following fragmentation process \cite{DeanMaj}.
Fix $b\ge2$  and a random vector $\vv=\vvb$.
%This random vector will be fixed in this paper.
Note that the definitions and results below depend only on the
distribution of \vvb, so it would be more precise to say that we fix a
distribution on $\bbR^b$; we find it, however, more convenient to
state the results in terms of a random vector.
We assume throughout the paper that
$0\le V_j\le1$, $j=1,\dots,b$, and
\begin{equation}
  \label{a1}
\sumjb V_j=1,
\end{equation}
\ie, that \vvb{} belongs to the standard simplex.
For simplicity we also assume that each $V_j<1$ a.s.
We allow $V_j=0$, but note that, a.s., $0<V_j<1$ for at least one $j$.

Starting with an object of size $x\ge1$, we break it into $b$ pieces
with sizes $V_1x,\dots,V_bx$.
Continue recursively with each piece of size $\ge1$, using new
(independent) copies of the random vector \vvb{} each time.
The process terminates a.s.\ after a finite number of steps,
leaving a finite set of fragments of sizes $<1$.
We let $\nx$ be the random number of fragmentation events, i.e., the
number of pieces of size $\ge1$ that appear during the process;
further, let $\nex$ be the final number of fragments, i.e., the number
of pieces of size $<1$ that appear.

This model has been studied by Dean and Majumdar \cite{DeanMaj}, who
found (without giving a rigorous proof) that the asymptotic behaviour
of $\nx$ as $\xtoo$ depends on the position of the roots of a certain
characteristic equation; we give a precise version of this in
\refT{T1} below. Some special cases have earlier been studied by other
authors, see \refS{Sex}.

It is natural to consider the fragmentation process as a tree, with
the root representing the original object, its children the result of
the first fragmentation, and so on.
It is then convenient to let the fragmentation go on for ever,
although we ignore what happens to pieces smaller than 1.
Let us label each node with the size of the corresponding object.

We thus consider the infinite rooted $b$-ary tree
$T_b$, whose nodes are the strings $J=\jj{k}$ with
$j_i\in\set{1,\dots,b}$ and $k\ge0$. Let $\cBx$
denote the set of all such strings, and let
$(V_1^{(J)},\dots,V_b^{(J)})$, $J\in\cBx$, be
independent copies of \vv. Then node $J=\jj{k}$
gets the label $x\prod_{i=1}^k
V_{j_i}^{(\jj{i-1})}$. Thus \nx{} is the number of
nodes with labels $\ge1$, \ie
\begin{equation}
  \label{a4}
\nx
=\sum_{J\in\cBx}\Bigind{
{x\prod_{i=1}^k V_{j_i}^{(\jj{i-1})}\ge1}}.
\end{equation}

By the recursive construction of the fragmentation process,
we have $\nx=0$ for $0\le x<1$ and
\begin{equation}
  \label{a5}
\nx\eqd 1+\sumjb N^{(j)}(V_j x),
\qquad x\ge1,
\end{equation}
where $N^{(j)}(\cdot)$ are copies of the process $N(\cdot)$,
independent of each other and of \vvb.

\begin{remark}
Let $\tx$ be the subtree of $\too=T_b$ consisting of all nodes with
labels $\ge1$. Then $\nx=|\tx|$, the number of nodes in \tx. More
precisely, we call these nodes \emph{internal nodes} of \tx, and we
say that a node in \too{} is an \emph{external node} of \tx{} if it
has label $<1$ but
  its parent is an internal node.

Thus $\nx$ is the number of internal nodes, and
\nex{} is the number of external nodes. Since each internal node
has $b$ internal or external children, we have, for $x\ge1$,
$\nx+\nex=1+b\nx$, or $\nex=(b-1)\nx+1$.
Hence the results for \nx{} immediately yield similar results for
\nex{} and the total number of external and internal nodes $\nx+\nex$ too.

In this paper we thus study the size of the
fragmentation tree \tx. Of course, it is
interesting to study other properties too, such as
height, pathlength, \dots.

Note that we may define $\tx$ for all $x\ge0$ simultaneously, using the same
  $V_j^{(J)}$; this defines $(\tx)_{x\ge0}$ as an increasing stochastic
  process of trees. (Equivalently, we can label all nodes of $\too$
  as above, starting with 1 at the root, and then keep all nodes with
  labels $\ge1/x$).
\end{remark}

\begin{remark}
  We assume for convenience that each object is split into the same
  number $b$ of parts. Our method applies also to the case of a random
  number of parts. Indeed, if the number of parts is bounded, we can
  use the results below with $b$ large enough, setting the non-existing
$V_j\=0$.
If the number of parts is unbounded, we can, under suitable assumptions,
use the proofs below with minor modifications.
We leave this extension to the reader.
\end{remark}

Our main result is \refT{T1} below on the asymptotic distribution
of $\nx$, together with the
corresponding estimates for mean and variance given in \refT{TMV}.

We define (with $0^z\=0$), at least for $\Re z\ge0$,
\begin{equation}
  \label{phi}
\phi(z)\=\sumjb \E V_j^z,
\end{equation}
and note that $\phi(z)$ is bounded and analytic in the open
right half-plane \set{z:\Re z>0}.
Since we assume \eqref{a1}, clearly $\phi(1)=1$.
Since further $0\le V_j<1$ a.s., the function $\phi(z)$ is decreasing
for real $z>0$; hence $\phi(z)>1$ when $0<z<1$ and $\phi(z)<1$ for
$1<z<\infty$. Further, $|\phi(z)|\le\sum_j\E|V_j^z|=\phi(\Re z)$, so
$|\phi(z)|<1$ when $\Re z>1$.

A crucial role is played by the solutions to the \emph{\che}
\begin{equation}
  \label{chareq}
\phi(\gl)=1.
\end{equation}
By the comments above, $\gl=1$ is one root, and
$\Re\gl\le1$ for every root $\gl$; furthermore,
there is no real root in $(0,1)$.

We further define
\begin{equation}\label{alpha}
 \ga\=-\phi'(1)=\sumjb\E(-V_j\ln V_j),
\end{equation}
the expected entropy of \vvb.

We need a (weak) regularity condition on the distribution of \vvb. We
find the following convenient, although it can be weakened to
\conditionbd{} in \refS{Sprel} for suitable $\gd$.
For examples where this regularity and \refT{T1} fail, see \refE{Elattice}.
\begin{conditionA}
Each $V_j$ has a distribution that is absolutely continuous on
$(0,1)$, although a point mass at $0$ is allowed.
\end{conditionA}

Note that there is no condition on the joint
distribution. In one case, however, we need a
condition including the joint distribution too.
(Note that both conditions are satisfied if $\vv$
has a density on the standard simplex, \ie{} if
$(V_1,\dots,V_{b-1})$ has a density.)

\begin{conditionAAA}
The support of the distribution of $\vv$ on the standard simplex has an interior point.
\end{conditionAAA}

If \condition{A} holds, then, by Lemmas
\refand{LA}{LBd} below, there is only a finite
number of roots of $\phi(\gl)=1$ in
\set{\gl:\Re\gl\ge\gd} for any $\gd>0$. We may thus
order the roots with $\Re\gl>0$ as $\gl_1,
\gl_2,\dots,\gl_\Nx$ with decreasing real parts:
$\gl_1=1 > \Re\gl_2\ge\Re \gl_3\ge\dots$; we will
assume this in the sequel. If $\gl_1=1$ is the only
root with $\Re\gl>0$, we set $\gl_2=-\infty$ for
convenience.

We let ${\cM}^\mathbb{C}$ denote the space of
probability measures on $\bbC$, and let
\begin{align*}
 {\cM}^\mathbb{C}_2(\gamma)
 :=\{\eta\in {\cM}^\mathbb{C}\, : \,
 \int|z|^2\dd\eta(z)<\infty, \text{ and } \int z \dd\eta(z) =\gamma\},
% \|\eta\|_2<\infty, \text{ and } \E \eta =\gamma\},
 \qquad \gamma \in \mathbb{C}.
\end{align*} We let $T$ denote the map (assuming $\gl_2\neq-\infty$)
\begin{align}\label{limmap}
T: {\cM}^\mathbb{C}\to {\cM}^\mathbb{C},
 \quad \eta \mapsto {\cL}\left(\sum_{r=1}^b V_r^{\lambda_2} Z^{(r)}\right),
\end{align}
where $(V_1,\ldots,V_b)$, $Z^{(1)},\ldots,Z^{(b)}$
are independent and ${\cL}(Z^{(r)})=\eta$ for
$r=1,\ldots,b$. Note that $T$ maps
${\cM}^\mathbb{C}_2(\gamma)$ into itself for each
$\gamma$, since $\gl_2$ satisfies $\phi(\gl_2)=1$.

We state our main result. The constant $\ga>0$ is defined in
\eqref{alpha} above and $\gb$ is given explicitly in \refT{TMV}.
The $\ell_2$ distance between distributions is defined in \refS{Smetrics}.
\begin{theorem}
  \label{T1}
Suppose that \condition{A} holds.
Then we have:
\begin{thmxenumerate}
\item
If\/ $\Re\gl_2<1/2$ then $\E \nx = \gai x + o(\sqrt{x})$,
$\Var\nx\sim \gb x$ with $\gb>0$ and
\begin{align*}
\frac{\nx -\gai x}{\sqrt{ x}} \dto {\cN}(0,\gb).
\end{align*}
\item
If\/  $\Re\gl_2=1/2$ and each root $\gl_i$ with
$\Re\gl_i=1/2$ is a simple root of $\phi(\gl)=1$,
and further \condition{\AAA} too holds,  then $\E
\nx = \gai x + O(\sqrt{x})$, $\Var(\nx)\sim \gb
x\ln x$ with $\gb>0$ and
\begin{align*}
\frac{\nx -\gai x}{\sqrt{x\ln x}} \dto {\cN}(0,\gb).
\end{align*}
\item
If\/  $\Re\gl_2>1/2$, and $\gl_2$ and $\gl_3=\overline{\gl_2}$ are the only
roots of \eqref{chareq} with this real part, and these roots are simple,
then
$\E \nx = \gai x  + \Re(\gamma x^{\lambda_2})+O(x^\kappa)$,
with $\gai>0$, $\gamma \in \mathbb{C}\setminus\{0\}$, $1/2<\kappa<\Re\gl_2$ and
\begin{align*}
\ell_2\left(\frac{\nx  -\gai x}{x^{\Re\gl_2}},
\Re\bigpar{\Xi e^{\ii\Im\gl_2\ln x}}\right)
=O\left( x^{\kappa-\Re\gl_2}\right),
\end{align*}
for some complex random variable $\Xi$. Furthermore,
$\mathcal{L}(\Xi)$ is the unique fixed point of $T$ in
${\cM}^\mathbb{C}_2(\gamma)$.
\end{thmxenumerate}
\end{theorem}

\begin{remark}
We can regard our process as a general
(age-dependent) branching process \cite[Chapter
6]{Jagers}, provided we make a logarithmic change
of time as in \refS{Smean}. Indeed, there are two
versions. For internal nodes, the individuals in
the branching process live for ever, and give birth
at times $-\ln V_1,\dots,-\ln V_d$. For external
nodes, we have a splitting process where each
individual when it dies gives birth to new
particles with life lengths $-\ln V_1,\dots,-\ln
V_d$. For both versions, we obtain a super-critical
branching process with Malthusian parameter 1, but
the identity \eqref{a1} causes the asymptotics for
moments and distributions to be quite different
from typical super-critical branching processes.
\end{remark}

\begin{remark}
If $\gl_2$ in \eqref{limmap} is real, then the stable distributions
of index $1/\gl_2$ are fixed points of $T$. Note, however, that in our
case, $\gl_2$ never is real. Moreover, the fixed points we are
interested in have finite variance, and are thus quite different from
stable distributions.

For the related {Quicksort} fixed point equation, Fill and
Janson \cite{SJ134}
found a complete characterization of
the set of fixed points; in that case, all fixed points are formed by
combining certain stable distributions with the unique fixed point with mean 0
and finite variance.
\end{remark}

\begin{remark}
\condition{\AAA} is needed only in part (ii),
and is needed only to exclude the possibility that for each root $\gl_i$ with
$\Re\gl_i=1/2$, there is a complex constant $C_i$ such that
\begin{equation}\label{crazy}
\sumjb V_j^{\gl_i} = C_i \quad\text{a.s.}
\end{equation}
This is easily seen to be impossible if \condition{\AAA} holds, and even otherwise it
seems highly unlikely for any particular example, but it seems very likely that there are
examples satisfying \condition{A}
where $\vv$ is concentrated on a curve, say, such that \eqref{crazy} holds.
\end{remark}

We will prove the statements on mean and variance, with further
refinements,
in \refS{Smean}. To prove convergence in distribution, we will use a
continuous time version of the
contraction method.
We develop a general theorem, that we find to be of independent
interest, in \refS{Scontr}. This theorem is applied to our problem in
\refS{Sproof}.
Some examples are given in Sections \refand{Sex}{Snex}.

%\section*{Acknowledgment}
\begin{acks}
This research was initiated and largely done
during conferences in Oberwolfach, Frankfurt
and Vienna in August and September 2004.
We thank Luc Devroye, Jim Fill, Allan Gut and Hsien-Kuei Hwang for
valuable comments.
\end{acks}

\section{Further preliminaries}\label{Sprel}

We define (again with $0^z\=0$)
\begin{equation}
  \label{psi}
\psi(z,w)
\=\Cov\bbigpar{\sumjb  V_j^z,\sumjb  V_j^w}
%=\E\bbigpar{\bbigpar{\sumjb  V_j^z-\phi(z)}\bbigpar{\sumkb  V_k^w-\phi(w)}}
=\E\bbigpar{\sumjb  V_j^z\sumkb  V_k^w}-\phi(z)\phi(w)
.
\end{equation}
In particular,
$\psi(z,\bz)
=\E\bigl|\sumjb  V_j^z-\phi(z)\bigr|^2\ge0$,
with equality only if $\sumjb V_j^z=\phi(z)$ a.s.

For $\Re z,\Re w\ge0$, we have $|V_j^z|,|V_j^w|\le1$ and thus
%$|\phi(z)|\le b$ and, if also $\Re w\ge0$,
$|\psi(z,w)|\le 2b^2$.

We say that \vv{} is \emph{lattice} if there exists a number $r$ with
$0<r<1$ such that every $V_j\in\set{r^n}_{n\ge0}\cup\set0$ a.s.;
otherwise \vv{} is \emph{non-lattice}. Basic Fourier analysis applied
to the probability measure $\nu$ defined in \eqref{qnu} shows that \vv{} is
non-lattice if and only if $\gl=1$ is the only root of
\eqref{chareq} with $\Re\gl=1$.
(Otherwise, there is an infinite number of roots with $\Re\gl=1$.)
We will assume this, and more, below.

We introduce a family of regularity conditions that are weaker than
\condition{A}.

\begin{conditionBd} (Here $\gd$ is a real number with $\gd\ge0$.)
  \begin{equation*}
    \limsup_{t\to\infty} |\phi(\gd+\ii t)|<1.
  \end{equation*}
\end{conditionBd}

\begin{lemma}
  \label{LBd}
If \condition{B($\gd$)} holds for some $\gd\ge0$, then
\condition{B($\gd'$)} holds for every $\gd'>\gd$ as well; moreover
  \begin{equation*}
    \limsup_{\substack{\Re z\ge \gd\\\Im z\to\infty}} |\phi(z)|<1.
  \end{equation*}
\end{lemma}

\begin{proof}
  Choose first $\eps>0$ such that
$\limsup_{t\to\infty} |\phi(\gd+\ii t)|<1-2\eps$,
and then $A$ such that
$|\phi(\gd+\ii t)|\le 1-2\eps$ if $t\ge A$, and thus also if $t\le-A$.
Recall further that $|\phi(\gd+\ii t)|\le b$ for all $t$.
Since $\phi(z)$ is  analytic, and thus harmonic, in the half-plane
$\chd\=\set{z:\Re z>\gd}$ and bounded and continuous in
$\bchd$, $\phi$ is given by the Poisson integral of its boundary values
\cite[Lemma 3.4]{Garnett}:
\begin{equation}
  \label{poi}
\phi(x+\ii y)=\intooo P_{x-\gd}(y-t)\phi(\gd+\ii t)\dd t,
\qquad x>\gd,
\end{equation}
where $P_x(y)=x/\xpar{\pi(x^2+y^2)}$, the Poisson
kernel for the right half-plane.  Let $\go(x+\ii
y)\=\int_{-A}^A P_{x-\gd}(y-t)\dd t$, the harmonic
measure of $[\gd-\ii A,\gd+\ii A]$; then
\eqref{poi} implies
\begin{equation}\label{poi2}
|\phi(x+\ii y)|\le \intooo P_{x-\gd}(y-t)|\phi(\gd+\ii t)|\dd t
\le b\go(x+\ii y)+1-2\eps.
\end{equation}
It is well-known, and easy to see, that the set
$B\=\set{z\in\chd:\go(z)>\eps/b}$ is bounded; in fact, it is the
intersection of $\chd$ and a circular disc \cite[p.\ 13]{Garnett}.
Thus, $A_1\=\sup\set{\Im z:z\in B}<\infty$, and if $\Re x\ge\gd$ and
$|y|>A_1$, then $\go(z)\le\eps/b$ and
\eqref{poi2} yields $|\phi(x+\ii y)|\le 1-\eps$.
\end{proof}

\begin{lemma}\label{LA}
  If \condition{A} holds, then
\condition{B($\gd$)} holds for every $\gd\ge0$.
\end{lemma}

\begin{proof}
We have $\E V_j^{\ii t} = \E \bigpar{e^{\ii t\ln V_j}\ind{V_j>0}}$, the
Fourier transform of the distribution of $\ln V_j$
(ignoring any point mass at $0$), so by
\condition{A} and the Riemann--Lebesgue lemma, $\E V_j^{\ii t}\to0$ as \ttoo{}
for every $j$, and thus $\phi(\ii t)\to0$ as \ttoo.
Hence, \condition{B(0)} holds, and the result follows by \refL{LBd}.
\end{proof}

\begin{lemma}
  \label{bdfinite}
If \condition{B($\gd$)} holds for some $\gd>0$, then there is only a
finite number of roots to $\phi(\gl)=1$ with $\Re\gl\ge\gd$.
\end{lemma}

\begin{proof}
  By \refL{LBd}, all such roots satisfy $|\Im\gl|\le C$ for some
  $C<\infty$.
Furthermore, all roots satisfy $\Re\gl\le1$, so if further
  $\Re\gl\ge\gd$, $\gl$ belongs to a compact rectangle $K$ in the open
  right half-plane. Since, $\phi(z)-1$ is analytic and non-constant in
  this half-plane,
  it has only a finite number of roots in $K$.
\end{proof}

In particular, by the comments above,
\conditionbd{} with $\gd\le1$ implies that \vv{} is non-lattice.

\section{Mean and variance}\label{Smean}

We let $\gL$ denote the set of solutions to the \che{} \eqref{chareq},
\ie
\begin{equation}
\label{Lambda}
  \gL\=\set{\gl:\phi(\gl)=1};
\end{equation}
we further define its subsets
\begin{equation}
\label{Lambdas}
  \gL(s)\=\set{z\in\gL:\Re(z)=s}.
\end{equation}
In general, $\phi(\gl)$ is defined only for
$\Re\gl\ge0$,  and we consider only such $\gl$ in
\eqref{Lambda}. However, in cases where $\phi$
extends to a meromorphic function in a larger
domain (for example, when $\phi$ is rational), we
may include such $\gl$ too in $\gL$; this makes no
difference in \refT{TMV}. (In \refT{Trational}, we
include all roots in the complex plane.) We will
use $\gL(s)$ only for $s\ge0$, where there is no
ambiguity.

Let $m(x)\=\E\nx$ and $\gss(x)\=\Var\nx$.
We will show the following asymptotics.

\begin{theorem}
  \label{TMV}
Assume that \condition{B($\gd$)} holds with $0\le\gd<1$,
and let
$\gl_1,\dots,\gl_\Nx$ be the elements of \set{\gl\in\gL:\Re\gl>\gd},
ordered with $\gl_1=1>\Re \gl_2 \ge \Re \gl_3 \ge \cdots$.
Then, the following hold as \xtoo:
%with $\gsx\=\Re\gl_2$,
\begin{thmxenumerate}
  \item
$m(x)\sim \gai x$.
\item
If further $\phi'(\gl_i)\neq0$ for $i=1,\dots,\Nx$,
\ie{}, each $\gl_i$ is a simple root of $\phi(\gl)=1$,
then, more precisely,
for every $\gd'>\gd$,
\begin{equation}\label{tmv2}
m(x)=\sum_{i=1}^\Nx \frac{1}{-\gl_i\phi'(\gl_i)} x^{\gl_i} + O(x^{\gd'}).
\end{equation}
\item
If\/ $\gd<1/2$ and either  $\Nx=1$ or $\Re \gl_2<1/2$, then
$\gss(x)\sim \gb x$, with
\begin{equation}\label{tmv}
  \gb=\ga\qi\frac{1}{2\pi}\int_{-\infty}^\infty
 \frac{\psi(\hiu,\hiuu)}{|\hiu|^2|1-\phi(\hiu)|^2} \dd u
\in(0,\infty).
\end{equation}
\item
If\/ $\Nx\ge2$ and $\Re \gl_2=1/2$,
and each $\gl_i$ with $\Re\gl_i=1/2$ is a simple root of $\phi(\gl)=1$,
then
$\gss(x)= \gb x\ln x+o(x\ln x)$, with
\begin{equation} \label{julie}
  \gb=\sum_{\gl\in\gL(1/2)}
\frac1{\ga|\gl\phi'(\gl)|^2} \psi(\gl,\overline{\gl})
\ge0.
\end{equation}
If, moreover, \condition{\AAA} holds  (or, more
generally, for some $\gl_i\in\gL(1/2)$,
\eqref{crazy} does not hold), then $\gb>0$.
\item
If\/ $\Nx\ge2$ and $\Re \gl_2>1/2$,
and each $\gl_i$ with $\Re\gl_i=\Re\gl_2$ is a simple root of $\phi(\gl)=1$,
then
\begin{equation*}
  \gss(x)=\!\!\!\!\sum_{\gl_i,\gl_k\in\gL(\Re\gl_2)}
\frac1{\gl_i\gl_k\phi'(\gl_i)\phi'(\gl_k)\bigpar{1-\phi(\gl_i+\gl_k)}}
\psi(\gl_i,\gl_k) x^{\gl_i+\gl_k} + o\bigpar{x^{2\Re\gl_2}}.
\end{equation*}
\end{thmxenumerate}
\end{theorem}

\begin{remark}\label{Rnl}
It follows from the proof that for (i) we do not need \condition{B($\gd$)};
it is enough that \vv{} is non-lattice.
\end{remark}

\begin{remark}
The case when some $\phi'(\gl_i)=0$ is similar; now terms
$x^{\gl_i}\ln x$ (and possibly $x^{\gl_i}\ln^d x$, $d\ge2$) will
appear in \eqref{tmv2}. We leave the details to the reader.
\end{remark}

If $\phi$ is a rational function,  then \eqref{tmv2} can be improved
to an exact formula.
Furthermore, in case (iii) of \refT{TMV} we then can give an alternative
formula for $\gb$.

\begin{theorem}
  \label{Trational}
Assume that $\phi$ is (\ie, extends to) a rational fumction, and let
$\gl_1,\dots,\gl_\Nx$ be the roots of $\phi(\gl)=1$ in the complex
plane. Suppose further that all these roots are simple.
  \begin{thmxenumerate}
\item Then
\begin{equation}\label{tmv2x}
m(x)=\sum_{i=1}^\Nx \frac{1}{-\gl_i\phi'(\gl_i)} x^{\gl_i}
-\frac1{\phi(0)-1},
\qquad x\ge1.
\end{equation}
\item
Assume further that $\gl_1=1$ and $\Re\gl_i<1/2$ for $i=2,\dots,\Nx$,
and that $V_j>0$ a.s.\ for every $j$.
Define, for notational convenience, $\gl_0\=0$, $a_0\=-1/(b-1)$ and
$a_i\=-1/(\gl_i\phi'(\gl_i))$ for $i=1,\dots,\Nx$.
Then $\gss(x)\sim \gb x$, with
\begin{align*}
  \gb
&=
\gai\sum_{i,k\neq1}\frac{a_ia_k}{1-\gl_i-\gl_k}
\Bigpar{\sum_{j,l=1}^b \E V_j^{\gl_i}V_l^{\gl_k} (V_j\bmin V_l)^{1-\gl_i-\gl_k}
%-2\sum_{j=1}^b\E V_j^{1-\gl_k} +1
-2\phi(1-\gl_k) +1
}
\\&\qquad
-2\ga^{-2}\sum_{i=2}^\Nx\frac{a_i}{\gl_i}
\Bigpar{\sum_{j,l=1}^b\E (V_j^{\gl_i}V_l^{1-\gl_i}-V_l)\ind{V_l\le V_j}
-\phi(1-\gl_i) +1
}
\\&\qquad
-2\ga^{-2}a_0
\Bigpar{\sum_{j,l=1}^b
\E V_l(\ln V_j-\ln V_l)\ind{V_l< V_j}-\ga }
\\&\qquad
+\ga^{-3}
\Bigpar{\sum_{j,l=1}^b
\E (V_j\bmin V_l)-1
}
-\gai.
\end{align*}
  \end{thmxenumerate}
\end{theorem}

The proof of these theorems will occupy the remainder of this section.
We first show that all  moments of $\nx$ are finite.
\begin{lemma}
  \label{Lmom}
For every $m\ge1$ and $x\ge0$, $\E\nx^m<\infty$.
Furthermore, $\sup_{0\le y\le x}\E\ny^m<\infty$.
\end{lemma}
\begin{proof}  For a string  $J=j_1\cdots j_k\in \cBx$ we denote by
  $|J|=k$ the depth of the corresponding node in $T_b$. Note that we
  have $|\{J\in\cBx:0\le |J|\le k\}|\le b^{k+1}$. Hence, if
  $N(x)>b^{k+1}$ for some $k\ge 0$ then by \eqref{a4} there exists a $J=j_1\cdots
  j_k\in \cBx$
with $x\prod_{i=1}^k V_{j_i}^{(\jj{i-1})}\ge 1$.
Markov's inequality implies that for all $q\ge 1$
\begin{align*}
\P (N(x)>b^{k+1}) &\le \P \left(\bigcup_{J\in
\cBx:|J|=k}\left\{
    \prod_{i=1}^k V_{j_i}^{(\jj{i-1})}\ge 1/x\right\}\right)\\
& \le \sum_{J\in \cBx:|J|=k}
\P\left(   \prod_{i=1}^k V_{j_i}^{(\jj{i-1})}\ge 1/x\right) \\
&\le \sum_{J\in \cBx:|J|=k} x^q \E  \prod_{i=1}^k
\left(V_{j_i}^{(\jj{i-1})}\right)^q\\
&= x^q \phi(q)^k.
\end{align*}
Hence, for all $y\ge b$ we obtain with  $k=\lfloor
\log_b y\rfloor-1$ and $\phi(q)\le 1$ that
\begin{equation}\label{nnn1}
\P (N(x)>y) \le  x^q \phi(q)^k \le x^q
\phi(q)^{\log_b y -2}=\frac{x^q}{\phi(q)^2}\,
y^{\log_b \phi(q)}.
\end{equation}
We have $\phi(q)\to 0$ as $q\to\infty$ since
$V_j<1$ a.s.~and  by dominated convergence. Hence,
for all  $m\ge 1$ there  exists a $q>0$ with
$\log_b \phi(q)<-m$.  The tail bound \eqref{nnn1}
thus implies $\E N(x)^m<\infty$ for all $m\ge 1$
and all $x\ge 0$.

The final statement follows because $0\le
\ny\le\nx$ when $0\le y\le x$.
\end{proof}

We find it convenient to switch from multiplicative
to additive notion. We therefore define
\begin{align*}
  X_j&\=-\ln V_j\in(0,\infty], \qquad j=1,\dots,b,
\\
  \nxt&\= N(e^t), \qquad -\infty\le t<\infty.
\end{align*}

The definition \eqref{a4} and the recursive equation \eqref{a5}
thus translate to
\begin{align}
  \label{a4x}
\nxt
&
=\sum_{J\in\cBx}\Bigind{
{\sum_{i=1}^k X_{j_i}^{(\jj{i-1})}\le t}},
\\
  \label{b1}
\nxt
&\eqd 1+\sumjb \nxx^{(j)}(t-X_j),
\qquad t\ge0,
\end{align}
where $\nxx^{(j)}(\cdot)$ are independent copies of the process $\nxx(\cdot)$,
and $\nxt=0$ for $-\infty\le t<0$.
Further define
\begin{align*}
  \mxt&\=\E\nxt=m(e^t),
\\
  \gssxt&\=\Var\nxt=\gss(e^t).
\end{align*}
Thus $\mxt=\gssxt=0$ for $t<0$. Taking expectations in \eqref{b1} we
find
\begin{equation}
  \label{b2w}
\mxt=1+\E\sumjb \mxx(t-X_j),
\qquad t\ge0.
\end{equation}

Let $\mu_j$ be the distribution of $X_j$ on $(0,\infty)$;
this is a measure of mass $1-\P(V_j=0)$;
let further $\mu\=\sumjb\mu_j$.
Then \eqref{b2w} can be written
\begin{equation}
  \label{sofie}
\mxt=1+\sumjb\mu_j*\mxt = 1+\mu*\mxt,
\qquad t\ge0,
\end{equation}
where $\mu*f(t)=\int_0^\infty f(t-x)\dd\mu(x)$.
This is the standard renewal equation, except that $\mu$ is not a
probability measure.

Similarly, conditioning on $X_1,\dots,X_b$, for $t\ge0$,
\begin{equation*}
  \begin{split}
\hskip6em&\hskip-6em
\E\bigpar{(\nxt-\mxt)^2\mid X_1,\dots,X_b}
\\
&=
\E\Bigpar{\Bigparq{\sumjb \nxx^{(j)}(t-X_j)+1-\mxt}^2\Bigm|
X_1,\dots,X_b}
\\
&=
\Var\Bigpar{\sumjb \nxx^{(j)}(t-X_j)+1-\mxt\Bigm|X_1,\dots,X_b}
\\
&\hskip10em
+ \Bigpar{\sumjb\mxx(t-X_j)+1-\mxt}^2
\\
%=
%\E\Bigpar{\Bigparq{\sumjb \bigpar{\nxx^{(j)}(t-X_j)-\mxx(t-X_j)}
%+1+\sumjb\mxx(t-X_j)-\mxt}^2\Bigm| X_1,\dots,X_b}
%\\
&=
\sumjb\gssx(t-X_j)+\Bigpar{\sumjb\mxx(t-X_j)-\mxt+1}^2.
  \end{split}
\end{equation*}
Taking the expectation we obtain
\begin{equation}
\label{emma}
  \gssxt
=
\E\sumjb\gssx(t-X_j)+h(t)
=\mu*\gssx(t)+h(t),
\qquad t\ge0,
\end{equation}
where, recalling \eqref{b2w},
\begin{align}
  h(t)&\=\E\Bigpar{\sumjb\mxx(t-X_j)-\mxt+1}^2
\notag
\\
&=\E\Bigpar{\sumjb\mxx(t-X_j)-\mxt}^2
+2\Bigpar{\E\sumjb\mxx(t-X_j)-\mxt}+1
\notag
\\&
=
\E\Bigpar{\sumjb\mxx(t-X_j)-\mxt}^2-1
.%,\qquad t\ge0.
\label{samuel}
\end{align}

Both \eqref{sofie} and \eqref{emma} are instances of the general
renewal equation \eqref{magnus} below, and from renewal theory we get
the following result. We say that a function on $\ooo$ is
\emph{locally bounded} if it is bounded on every finite interval.

\begin{lemma}
  \label{L1}
Assume that \vv{} is non-lattice.
Let $f$ be a locally bounded measurable function on $\ooo$.
Then the renewal equation
\begin{equation}
  \label{magnus}
F=f+\mu*F
\end{equation}
has a unique locally bounded solution $F$ on \ooo{}.
If further \vv{} is non-lattice, then
we have the following asymptotical results, as \ttoo,
\begin{thmxenumerate}
  \item
If\/ $f$ is a.e.\ continuous and
$\int_0^\infty \ff(t)\dd t<\infty$, where
$\ff(t)\=\sup_{u\ge t}e^{-u}|f(u)|$,
then
$F(t)=(\gamma+o(1))e^t$,
with $\gamma=\ga\qi\intoo f(t)e^{-t}\dd t$.
\item
If\/ $f(t)=e^{t}$, then $F(t)\sim \gai t e^t$.
\item
If\/ $f(t)=e^{\gl t}$ with\/ $\Re\gl=1$ and\/ $\Im\gl\neq0$, then
$F(t)=o\bigpar{t e^{t}}$.
\item
If\/ $f(t)=e^{\gl t}$ with $\Re\gl>1$, then
$F(t)\sim\bigpar{1-\phi(\gl)}\qi e^{\gl t}$.
\end{thmxenumerate}
\end{lemma}

\begin{proof}
For a function $f$ on $(0,\infty)$ and $z\in\bbC$, we define, when the integral
exists, the Laplace transform
$\tf(z)\=\intoo e^{-zt}f(t)\dd t$.
Similarly, the Laplace transform of $\mu$ is
\begin{equation}
\label{b2x}
  \tmu(z)\=\intoo e^{-tz}\dd\mu(t)
=
\sumjb\E e^{-zX_j}
=\sumjb\E V_j^z
=\phi(z),
\end{equation}
at least for $\Re z\ge0$.

 Since $\mu$ is not a probability measure, we define another
 (``conjugate'' or ``tilted'') measure $\nu$ on $[0,\infty)$ by
   \begin{equation}
\label{qnu}
\dd\nu(x)=e^{-x}\dd\mu(x).
   \end{equation}
Then $\nu$ is a probability measure because, by \eqref{a1},
\begin{equation*}
\nu[0,\infty)
=\intoo e^{-x}\dd\mu(x)
=\sumjb\intoo e^{-x}\dd\mu_j(x)
=\sumjb\E e^{-X_j}
=\sumjb\E V_j
=1.
\end{equation*}
Further, the mean of the distribution $\nu$ is
\begin{equation}\label{b2b}
\E  \nu
=\intoo x\dd\nu(x)
=\intoo xe^{-x}\dd\mu(x)
%=\sumjb\intoo e^{-x}\dd\mu_j(x)
=\sumjb\E \bigpar{X_je^{-X_j}}
=\sumjb\E \bigpar{(-\ln V_j)V_j}
=\ga
\end{equation}
and the Laplace transform is, for $\Re z\ge0$,
recalling \eqref{b2x},
\begin{equation}\label{b2c}
\tnu(z)
\=\intoo e^{-zx}\dd\nu(x)
=\intoo e^{-x-zx}\dd\mu(x)
=\tmu(z+1)
=\phi(z+1).
\end{equation}

Let $g(t)\=e^{-t}f(t)$ and $G(t)\=e^{-t}F(t)$. Then \eqref{magnus}
translates to
\begin{equation*}
  \begin{split}
  G(t)
&=e^{-t}F(t)
=e^{-t}f(t)+\intoo e^{-t}F(t-x)\dd\mu(x)
\\&
=g(t)+\intoo G(t-x)e^{-x}\dd\mu(x)
=g(t)+\nu*G(t).
  \end{split}
\end{equation*}
In other words, $G$ satisfies the renewal equation for the probability
measure $\nu$, so we can use standard results from renewal theory.

First, it is well known that the equation $G=g+\nu*G$ has a unique
locally bounded solution
which is given by $G=\sumn\nu^{*n}*g$, and thus
$F=\sumn\mu^{*n}*f$;
see \eg{} \cite[Theorem IV.2.4]{Asm} (which also applies directly to
$F$).
If we let $Y_1,Y_2,\dots$ be \iid{} random variables with the
distribution $\nu$, and let $S_n\=\sum_1^n Y_i$, this can be written
\begin{equation}\label{gg}
  G(t)=\sumn\E\bigpar{g(t-S_n)\ind{S_n\le t}}
=\E\sum_{S_n\le t} g(t-S_n).
\end{equation}

Under the assumptions of (i), $\ff$ is
non-increasing and integrable;
further, $\sup\ff\le \sup_{[0,1]}|f|+\ff(1)<\infty$, so $\ff$ is bounded
too. Hence \cite[Proposition IV.4.1(v),(iv)]{Asm} shows that
$\ff$ and $g$ are directly Riemann integrable.
The key renewal theorem \cite[Theorem IV.4.3]{Asm} and \eqref{b2b} now yield
$G(t)\to\gai\intoo g(x)\dd x=\gamma$, which proves (i).

In case (ii) we have $g(t)=1$, and thus $G(t)\sim\gai t$ by the
elementary renewal theorem \cite[IV.(1.5) and Theorem 2.4]{Asm}.

For (iii), $g(t)=e^{(\gl-1)t}=e^{ibt}$ for some real $b\neq0$.
The solution to \eqref{magnus} may be written
\cite[Theorem IV.2.4]{Asm}
%\cite[Theorem VI.6.1]{Feller}
$G(t)=\int_0^tg(t-x)\dd U(x)$, where $U$ is the locally bounded
solution to $U=1+\nu*U$ (\ie, $U=G$ for case (ii)).
Since, in analogy with \eqref{b2b}, $\int
x^2\dd\nu(x)=\sum_j\E\bigpar{(\ln V_j)^2V_j}<\infty$, the distribution
$\nu$ has finite variance, and the renewal theorem has the sharper
version \cite[Theorem XI.3.1]{Feller}
\begin{equation*}
  U(t)=\ga\qi t+c+R(t),
\end{equation*}
where $c$ is a certain constant ($\int x^2\dd\nu/2\ga^2$) and
$R(t)\to0$ as \ttoo.
Hence, using integration by parts for one term,
\begin{equation*}
  \begin{split}
  G(t)
&=\int_0^te^{\ii b(t-x)} \ga\qi\dd x+ce^{\ii bt}
+\int_0^te^{\ii b(t-x)} \dd R(x)
\\&
= O(1)+R(t)-R(0)e^{\ii bt} +\ii b\int_0^t e^{\ii b(t-x)} R(x)\dd x
=o(t).
  \end{split}
\end{equation*}

For (iv), we have by \eqref{gg} with $g(t)=e^{(\gl-1)t}$,
using dominated convergence and \eqref{b2c},
\begin{align*}
  e^{-\gl t} F(t)
&
=
e^{(1-\gl) t} G(t)
=
\sum_{n=0}^\infty\E \Bigpar{e^{(1-\gl)t+(\gl-1)(t-S_n)}\ind{t\le S_n}}
\\&
\to
\sum_{n=0}^\infty\E \bigpar{e^{-(\gl-1)S_n}}
=
\sum_{n=0}^\infty\E \bigpar{e^{-(\gl-1)Y_1}}^n
=
\sum_{n=0}^\infty\tnu(\gl-1)^n
\\&
=
\sum_{n=0}^\infty\phi(\gl)^n
=\bigpar{1-\phi(\gl)}\qi.
\end{align*}
\vskip-\baselineskip
\end{proof}

\begin{proof}[Proof of \refT{TMV}]
We first apply \refL{L1}(i) to \eqref{sofie}, with
$f(t)=1$ for $t\ge0$,
and obtain $\gamma=\gai$ and $\mxt\sim\ga\qi e^t$, which proves \refT{TMV}(i).

To obtain more refined asymptotics, we use Laplace transforms.
Let $H(t)\=\ind{t\ge0}$ (the Heaviside function),
and note that $\tH(z)=\intoo e^{-tz}\dd t=1/z$, $\Re z>0$.
Since the Laplace transform converts convolutions to products, the
%renewal equation \eqref{magnus} yields
%$\tF(z)=\tf(z)+\tmu(z)\tF(z)$, and thus
%\begin{equation}
%  \tF(z)=\frac{\tf(z)}{1-\tmu(z)},
%\end{equation}
%for $z$ such that the transforms exist.
renewal equation \eqref{sofie} yields
$\tmxx(z)=\tH(z)+\tmu(z)\tmxx(z)$, and thus
\begin{equation}\label{erika}
  \tmxx(z)=\frac{\tH(z)}{1-\tmu(z)} =\frac{1}{z(1-\phi(z))},
\end{equation}
for $z$ such that the transforms exist.
By the estimate $\mxt\sim\gai e^t$ above,
$\mxt=O(e^t)$ and thus $\tmxx(z)$ exists for $\Re z>1$.
Consequently, \eqref{erika} holds for
$\Re z>1$, and can be used to extend $\tmxx(z)$ to a meromorphic
function for $\Re z>0$.

We want to invert the Laplace transform in \eqref{erika}.
This is simple if $\phi$ is rational, yielding
\eqref{tmv2x}. (Note that $\phi(0)=\E|\set{j:V_j>0}|>1$.)
In general,
there are difficulties to doing this directly, because
$\tmxx(z)$ is not integrable along a vertical line $\Re z=s$; it
decreases too slowly as $|\Im z|\to\infty$.
We therefore regularize. Let $\eps>0$, and let
$\he\=H*\eps\qi\ett_{[0,\eps]}$; thus
\begin{equation*}
  H_\eps(t)=
  \begin{cases}
0,& t<0, \\
1-t/\eps, & 0\le t<\eps, \\
1, & t\ge\eps.
  \end{cases}
\end{equation*}
%$H_\eps(t)=0$ for $t<0$, $1-t/\eps$ for $0\le t<\eps$, and $1$ for $t\ge\eps$.
Let $\mex=\sumn\mu^{*n}*\he$ be the locally bounded solution to
$\mex=\he+\mu*\mex$. Note that $\he(t)\le H(t)\le\he(t+\eps)$, and
thus
\begin{equation}
  \label{c1}
\mex(t)\le\mxt\le\mex(t+\eps).
\end{equation}
We have
\begin{equation*}
  \tHe(z)=\tH(z)\eps\qi\int_0^\eps e^{-zt}\dd t
= \frac{1-e^{-\eps z}}{\eps z^2},
\qquad \Re z>0,
\end{equation*}
and we find, arguing as for \eqref{erika} above,
\begin{equation*}%\label{erika}
\tmex(z)=\frac{\tHe(z)}{1-\tmu(z)}
=\frac{1-e^{-\eps z}}{\eps z^2(1-\phi(z))},
\end{equation*}
first for $\Re z>1$, and then for $\Re z>0$, extending $\tmex$ to a
meromorphic function in this domain. This function decreases (using
\condition{B($\gd$)} and \refL{LBd})
as $|\Im z|^{-2}$ on vertical lines $\Re z=s\ge \gd$,
and is thus integrable there. Hence, the Laplace inversion formula
(a Fourier inversion) shows that for any $s>1$ and $t\ge0$,
\begin{equation}\label{cc}
  \mex(t)
=\frac1{2\pi\ii}\int_{s-\ii\infty}^{s+\ii\infty} e^{tz}\tmex(z)\dd z.
\end{equation}
We may, increasing $\gd$ a little if necessary, assume that
$\phi(z)=1$ has no roots with $\Re z=\gd$;
in cases (iii), (iv) and (v) we may similarly assume that each
$\gl\in\gL$ with $\Re\gl>\gd$ has $\phi'(\gl)\neq0$.
It is then easy to show,
using \condition{B($\gd$)} and \refL{LBd},
that we may shift the line of integration
in \eqref{cc} to $\Re z=\gd$ and obtain, for $0<\eps\le1$,
\begin{align*}
\mex(t)
&
=\frac1{2\pi\ii}\int_{\gd-\ii\infty}^{\gd+\ii\infty} e^{tz}\tmex(z)\dd z
+\sum_{i=1}^\Nx \Res_{z=\gl_i}\Bigpar{e^{tz}\tmex(z)}
\\&
=O\bbigpar{e^{t\gd}\int_{\gd-\ii\infty}^{\gd+\ii\infty}
\min\Bigpar{\frac1{|z|},\frac\eps{|z|^2}} \,|\!\dd z|
}
+\sum_{i=1}^\Nx \frac{e^{t\gl_i}}{-\gl_i\phi'(\gl_i)}
  \frac{1-e^{-\eps\gl_i}}{\eps\gl_i}
\\&
=\sum_{i=1}^\Nx \frac{e^{t\gl_i}}{-\gl_i\phi'(\gl_i)}\bigpar{1+O(\eps)}
+O\Bigpar{e^{t\gd}\bigpar{1+\ln\epsix}}.
\end{align*}
Now choose $\eps\=e^{-t}$; we then obtain
\begin{equation*}
\mex(t)
=\sum_{i=1}^\Nx \frac{e^{t\gl_i}}{-\gl_i\phi'(\gl_i)}
+O(1)
+O\Bigpar{e^{t\gd}\xpar{1+t}},
\qquad t\ge0.
\end{equation*}
Replacing $t$ by $t+\eps$, we obtain the same estimate for
$\mex(t+\eps)$, and thus \eqref{c1} yields
\begin{equation}\label{jesper}
\mxt
=\sum_{i=1}^\Nx \frac{e^{t\gl_i}}{-\gl_i\phi'(\gl_i)}
+O\bigpar{e^{t\gd'}},
\qquad t\ge0,
\end{equation}
which yields \refT{TMV}(ii).

For the estimates of the variance, we use \refL{L1} and \eqref{emma}.
It is easily seen (by dominated convergence)
that $h$ in \eqref{samuel} is a.e.\ continuous.
Choose $\gd'>\gd$ with $\gd<\gd'<\Re\gl_\Nx$; in case (iii) with $\Nx=1$, let
further $\gd'<1/2$.

Note that then \eqref{jesper} trivially holds for $t<0$ too.
%(provided $\gd'\le \Re\gl_\Nx$).
Hence,
\begin{align*}
\sumjb\mxx(t-X_j)-\mxt
&
=
\sum_{i=1}^\Nx \frac{\sumjb e^{(t-X_j)\gl_i}-e^{t\gl_i}}{-\gl_i\phi'(\gl_i)}
+O\bigpar{e^{t\gd'}}
\\&
=
\sum_{i=2}^\Nx \frac{e^{t\gl_i}}{-\gl_i\phi'(\gl_i)}
\bbigpar{\sumjb V_j^{\gl_i} -1}
+O\bigpar{e^{t\gd'}},
  \end{align*}
where we  use the fact that $\gl_1=1$ and thus
$\sumjb V_j^{\gl_1}-1=\sumjb V_j-1=0$.
Consequently, by \eqref{samuel} and \eqref{psi},
letting $\gsx\=\Re\gl_2>\gd'$  if $\Nx\ge2$,
and $\gsx\=\gd'$ if $\Nx=1$,
\begin{equation}\label{david}
  h(t)=
\sum_{i=2}^\Nx \sum_{k=2}^\Nx
\frac{e^{(\gl_i+\gl_k)t}}{\gl_i\gl_k\phi'(\gl_i)\phi'(\gl_k)}
\psi(\gl_i,\gl_k)
%\E\Bigpar{\sumjb V_j^{\gl_i} -1}\Bigpar{\sumjb V_j^{\gl_k} -1}
+O\bigpar{e^{t(\gd'+\gsx)}}.
\end{equation}

For \refT{TMV}(iii), \eqref{david} yields $h(t)=O\bigpar{e^{2\gsx t}}$
with $\gsx<1/2$,
and \refL{L1}(i) applies to \eqref{emma}, yielding $\gssx(t)\sim \gamma e^t$.
We postpone the calculation of $\gb=\gamma$, verifying \eqref{tmv}, to
\refL{L3}.

In \refT{TMV}(iv) and (v), we treat the terms in \eqref{david}
separately, using linearity; for the error term we also use
monotonicity and comparison with the case $f(t)=e^{t(\gd'+\gs_2)}$.
The results follow by applying \refL{L1}(i)--(iv) to \eqref{emma},
letting $t\=\ln x$,
with the leading terms coming from
the case $\Re\gl_k=\Re\gl_i=\Re\gl_2$, and for (iv) further
$\gl_k=\overline{\gl_i}$.

Furthermore, by \eqref{julie}, $\gb=0$ in (iv) only
if for  every $\gl_i\in\gL(1/2)$, we have
$\psi(\gl_i,\bar \gl_i)=\E|\sum_j
V_j^{\gl_i}-\phi(\gl_i)|^2=0$, and thus
\eqref{crazy} holds.
\end{proof}

\begin{lemma}
  \label{L3}
If\/ $\Re\gl_2<1/2$ or $\Nx=1$, then, with $h(t)$ as in \eqref{samuel},
\begin{equation*}
\intoo h(t)e^{-t}\dd t
=\frac{1}{2\pi}\int_{-\infty}^\infty
 \frac{\psi(\hiu,\hiuu)}{|\hiu|^2|1-\phi(\hiu)|^2} \dd u
>0.
\end{equation*}
\end{lemma}

\begin{proof}
Write $f(t)\=\bigpar{\mxx(t)-\gai e^t}e^{-t/2}$, $-\infty<t<\infty$.
Thus, by \eqref{jesper},
$f(t)=O\bigpar{e^{-(1/2-\gs_2)t}}$ for $t\ge0$
and
$f(t)=-\gai e^{t/2}=O\bigpar{e^{-|t|/2}}$ for $t<0$.
In particular, $f\in L^2(-\infty,\infty)$. Furthermore, the (two-sided)
Laplace transform
$\tf(z)\=\intooo f(t) e^{-tz}\dd t$ is analytic for $-(1/2-\gs_2)<\Re
z <1/2$.

Define further $f_1(t)\=f(t) e^{t/2} = \mxx(t)-\gai e^t$ and
$f_2(t)\=f_1(t)\ind{t\ge0}$.
Then $f_2(t)=O\bigpar{e^{\gs_2 t}}$, and thus the Laplace transform
$\tf_2(z)$ is analytic for $\Re z>\gs_2$. For $\Re z>1$ we have, by
\eqref{erika},
\begin{equation*}
  \tf_2(z)=\intoo e^{-tz}\bigpar{\mxx(t)-\gai e^t}\dd t
=\tmxx(z)-\gai(z-1)\qi
=\frac1{z(1-\phi(z)} -\frac1{\ga(z-1)};
\end{equation*}
by analytic continuation, this formula holds for $\Re z>\gs_2$.
Consequently, for $\gs_2<\Re z<1$,
\begin{equation*}
  \tf_1(z)=\tf_2(z)+\int_{-\infty}^0 e^{-tz} (-\gai e^t) \dd t
=\tf_2(z)-\gai(1-z)\qi
=\frac1{z(1-\phi(z)}.
\end{equation*}
Since $\tf(z)=\tf_1(z+1/2)$, we find the Fourier transform
\begin{equation}
  \label{manne}
\hatf(u)\=
\intoo e^{-\ii ut} f(t)\dd t
= \tf(\ii u)
= \tf_1(\tfrac12+\ii u)
=\frac1{(\tfrac12+\ii u)(1-\phi(\tfrac12+\ii u))} .
\end{equation}

Next, since $\sum_j e^{-X_j} = \sum_j V_j=1$,
\begin{equation*}
  \begin{split}
\sumjb \mxx(t- X_j) -\mxt
&= \sumjb f_1(t- X_j) -f_1(t) + \gai\sumjb e^{t-X_j} - \gai e^t
\\&
= \sumjb f_1(t- X_j) -f_1(t),
  \end{split}
\end{equation*}
so by \eqref{samuel}, and defining
$\Psi(x,y)\=\intooo f(t-x)f(t-y)\dd t$,
\begin{align*}
  \intoo h(t)e^{-t}\dd t+1
&
=  \intoo \bigpar{h(t)+1} e^{-t} \dd t
=
\E\intooo \bigpar{\sumjb f_1(t-X_j)-f_1(t)}^2 e^{-t} \dd t
\\&
= \E \intooo \Bigl| \sumjb e^{-X_j/2} f(t-X_j)-f(t)\Bigr|^2\dd t
\\&
= \E \sum_{j,k=1}^b e^{-X_j/2-X_k/2} \Psi(X_j,X_k)
  -2 \E \sumjb e^{-X_j/2} \Psi(X_j,0) + \Psi(0,0).
\end{align*}
By Parseval's relation and $f=\overline f$,
\begin{align*}
  \Psi(x,y)
=\frac1{2\pi} \intooo e^{-\ii ux} \hatf(u)\overline{e^{-\ii uy} \hatf(u)}
\dd u
\frac1{2\pi} \intooo |\hatf(u)|^2  e^{\ii u(y-x)} \dd u.
\end{align*}
Hence,
\begin{multline*}
  \intoo h(t)e^{-t}\dd t+1
=
\E\frac1{2\pi}\intooo |\hatf(u)|^2
\biggl(
 \sum_{j,k=1}^b e^{-X_j/2-X_k/2+\ii u(X_k-X_j)}
\\ \shoveright{
 - \sumjb e^{-X_j/2+\ii uX_j}   - \sumkb e^{-X_k/2-\ii uX_k}
  +1 \biggr) \dd u}
\\
\begin{aligned}
&=
\E\frac1{2\pi}\intooo |\hatf(u)|^2
  \E \Bigl|\sumjb V_j^{1/2-\ii u}-1\Bigr|^2 \dd u
\\&
=
\E\frac1{2\pi}\intooo |\hatf(u)|^2
 \bigpar{\psi(1/2+\ii u,1/2-\ii u)+|\phi(1/2+\ii u)-1|^2} \dd u.
\end{aligned}
\end{multline*}
Using \eqref{manne},
\begin{align*}
\E\frac1{2\pi}\intooo |\hatf(u)|^2 |\phi(1/2+\ii u)-1|^2 \dd u
=\intooo \frac{\dd u}{|1/2+\ii u|^2}
=\intooo \frac{\dd u}{\tfrac14+u^2}
=2\pi,
\end{align*}
and the formula follows by \eqref{manne}.

Since $\psi(z,\bar z)\ge0$, and by dominated convergence is continuous
for $\Re z\ge0$, it follows from \eqref{tmv} that $\gb=0$ only if for
every $z$ with $\Re z=1/2$,
$\psi(z,\bar z)=\E|\sum_j V_j^z-\phi(z)|^2=0$, and thus
\begin{equation}
  \label{gustaf}
\sumjb V_j^{1/2+\ii u} = \phi(1/2+\ii u)
\end{equation}
a.s., for every real $u$. Considering first rational $u$, we see that
\as{} \eqref{gustaf} holds for all real $u$.

However, for any realization $(V_1,\dots,V_b)$ and $\eps>0$, it is by
Dirichlet's theorem \cite[Theorem 201]{HW}
possible to find arbitrarily large $u$ with
$\Re V_j^{1/2+\ii u} \ge (1-\eps)V_j\qh$ for $j=1,\dots,b$,
and thus \eqref{gustaf} implies that
$\limsup_{u\to\infty} |\phi(\tfrac12+\ii u)| \ge \phi(1/2)\ge1$.
This contradicts \conditionbd{} and \refL{LBd}.
Hence $\gb>0$.
%Consider now a vector $(v_1,\dots,v_b)$ such that $v_j\ge0$ and
%$\sum_j v_j^{1/2+\ii u} = \phi(1/2+\ii u)$ for every $u\in\bbR$,
%\eg{} a.e. realization $(V_1,\dots,V_b)$.
%Let $J$ be a random index with the distribution
%$\P(J=j)=v_j\qh/\phi(1/2)$, and let $X$ be the random variable
%$\ln v_J$. Then
%\begin{align*}
%  \E e^{\ii u X}
%=\sumjb  e^{\ii u \ln v_j}\P(J=j)
%=\sumjb v_j^{1/2+\ii u}/\phi(1/2)
%= \phi(1/2+\ii u)/\phi(1/2).
%\end{align*}
%Hence $\phi$ determines the characteristic function of $X$, and thus
%its distribution; consequently, the set \set{v_1,\dots,v_b,0} is
%determined. Hence, \eqref{gustaf} implies that, a.s., each $V_j$ takes
%only at most $b+1$ values
\end{proof}

\begin{proof}[Proof of \refT{Trational}]
  \pfitem{i}
As remarked above, \eqref{tmv2x} follows by inverting the Laplace
transform in \eqref{erika}, using a partial fraction expansion.
\pfitem{ii}
Note first that $\phi(z)\to0$ as $z\to+\infty$ (by dominated
convergence), and thus, $\phi$ being rational, $\phi(\infty)=0$ and
$\phi(z)\to0$ as $|z|\to\infty$; hence \conditionbd{} holds for every
$\gd$.
We thus see that the conditions of \refT{TMV}(iii) are satisfied, and
from the proof above we see that, with $h$ given by \eqref{samuel},
\begin{equation*}
  \gb=\gai\intoo h(t)e^{-t}\dd t
=\gai \int_1^\infty h(\ln x) x^{-2}\dd x.
\end{equation*}
We have, by \eqref{tmv2x},
$m(x)=\sum_{i=0}^\Nx a_i x^{\gl_i}\ind{x\ge1}$ and thus
\begin{equation*}
  \mxx(\ln x-X_j)=m(xe^{-X_j})=m(xV_j)=\sum_{i=0}^\Nx
  a_i(xV_j)^{\gl_i}\ind{x\ge V_j\qi}.
\end{equation*}
Hence, letting $V_0\=1$, $\eps_0=-1$ and $\eps_j=1$ for $j\ge1$,
and recalling \eqref{a1},
\begin{align*}
  H(x)&\=
\sumjb \mxx(\ln x-X_j)-\mxx(\ln x)
=\sum_{i=0}^\Nx a_i x^{\gl_i}
\Bigpar{
\sumjb V_j^{\gl_i}\ind{x\ge V_j\qi}
-\ind{x\ge1}
}
\\&\phantom:
=\sum_{i\neq1} a_i x^{\gl_i} \sum_{j=0}^b  V_j^{\gl_i}\eps_j\ind{x\ge V_j\qi}
-a_1 x\sumjb  V_j\ind{1\le x< V_j\qi}.
\end{align*}
By \eqref{samuel}, this leads to
\begin{align*} \hskip4em&\hskip-2em
\int_1^\infty h(\ln x) x^{-2}\dd x +1
=
\int_1^\infty \E H(x)^2 x^{-2}\dd x
=
\E\int_1^\infty H(x)^2 x^{-2}\dd x
\\&
=
\E \sum_{i,k\neq1} a_ia_k \sum_{j,l=0}^b \eps_j\eps_l
V_j^{\gl_i}V_l^{\gl_k}
\int_{V_j\qi\bmax V_l\qi}^\infty x^{\gl_i+\gl_k-2}\dd x
\\&\qquad
-2 \E \sum_{i\neq1}a_1a_i\sum_{j=0}^b\sum_{l=1}^b\eps_j V_j^{\gl_i}V_l
\int_{V_j\qi}^{V_l\qi} x^{\gl_i-1}\dd x
\,\ind{V_j\qi\le V_l\qi}
\\&\qquad
+\E a_1^2\sum_{j,l=1}^b V_jV_l\int_1^{V_j\qi\bmin V_l\qi} \dd x
\end{align*}
and the result follows by straightforward calculations, noting that $a_1=\gai$.
\end{proof}

\section{Zolotarev metric and minimal $L_s$ metric}\label{Smetrics}

In this section we collect properties of the minimal $L_s$ metric
and the Zolotarev metric that are used subsequently.

We denote by $\mathcal{M}^d$ the space of probability measures on
$\Rset^d$.
 The \emph{minimal $L_s$ metric} $\ell_s$, $s>0$, is defined on the subspace
$\mathcal{M}^d_s\subset \mathcal{M}^d$ of probability measures with
finite absolute moment of order $s$ by
\begin{eqnarray*}
\ell_s(\mu,\nu):=\inf\left\{\|X-Y\|_s^{s\wedge1}\,:\,X\deq\mu,\,
Y\deq\nu\right\},\quad \mu,\nu \in \mathcal{M}^d_s,
\end{eqnarray*}
where $\|X\|_s\=(\E |X|^s)^{1/s}$
%where $\|X\|_s\=(\E |X|^s)^{(1/s) \wedge 1}$
denotes the $L_s$ norm of $X$. The infimum is taken
over all random vectors of $X$, $Y$ on a joint
probability space with the given marginal
distributions $\mu$ and $\nu$.  (In other words,
over all couplings $(X,Y)$ of $\mu$ and $\nu$.) We
will also use the notation
$\ell_s(X,Y):=\ell_s(\mathcal{L}(X),\mathcal{L}(Y))$.

For $s\ge 1$ and $\gamma\in \Rset^d$,  we denote by
$\mathcal{M}_s^d(\gamma)\subset \mathcal{M}^d_s$
the subspace of  probability measures with
expectation $\gamma$. The pairs
$(\mathcal{M}_s^d,\ell_s)$, $s>0$, and
$(\mathcal{M}_s^d(\gamma),\ell_s)$, $s\ge 1$, are
complete metric spaces and convergence in $\ell_s$
is equivalent to weak convergence plus convergence
of the absolute moments of order $s$.

Random vectors $(X,Y)$ with $X\deq \mu$, $Y\deq \nu$, and
$\ell_s(\mu,\nu)=\|X-Y\|_s^{s\wedge1}$ are called optimal couplings of
$(\mu,\nu)$. Such optimal couplings exist for all $\mu,\nu\in
\mathcal{M}^d_s$. These properties can be found in Dall'Aglio
\cite{DA56}, Major \cite{Ma78}, Bickel and Freedman \cite{BiFr}, and
Rachev \cite{rach2}.
Similar properties hold for probability measures
on $\mathbb{C}^d$ (because $\bbC^d\cong \bbR^{2d}$), where we use
corresponding notations.

The \emph{Zolotarev metric} $\zeta_s$, $s>0$ is defined by
\begin{equation}
\label{eq:3.6} \zeta_s(X,Y) \=
\zeta(\mathcal{L}(X),\mathcal{L}(Y))\= \sup_{f\in
\mathcal{F}_s}|\E(f(X) - f(Y))|
\end{equation}
where  $s=m+\alpha$ with $0<\alpha\le 1$,
$m=\ceil{s}-1\ge0$ is an integer,
and
\begin{eqnarray*}
\mathcal{F}_s:=\{f\in
C^m(\Rset^d,\Rset):\|f^{(m)}(x)-f^{(m)}(y)\|\le \|x-y\|^\alpha\},
\end{eqnarray*}
where $C^m(\Rset^d,\Rset)$ denotes the space of $m$ times
continuously differentiable functions $f$ on $\Rset^d$ and $f^{(m)}$
their $m$th derivative.

The expression $\zeta_s(X,Y)$ is finite if $X$ and $Y$ have finite
absolute moments of order $s$ and all mixed moments of orders
$1,\ldots,m$ of $X$ and $Y$ coincide.

The metric $\zeta_s$ is ideal of order $s$, i.e., we have for $Z$
independent of $(X,Y)$ and any $d\times d$ square matrix $A$
\begin{align*}
\zeta_s(X+Z,\,Y+Z)\le \zeta_s(X,Y), &&&
\zeta_s(AX,AY)\le\|A\|_{\op}^s\, \zeta_s(X,Y),
\end{align*}
where $\|A\|_\op\=\sup_{\|u\|=1}\|Au\|$ denotes the
operator norm of the matrix.  Convergence in
$\zeta_s$ implies weak convergence. For general
reference and properties of $\zeta_s$ we refer to
Zolotarev \cite{Zolotarev76,Zolotarev77} and Rachev
\cite{rach2}.

\section{General contraction theorems in continuous time}\label{Scontr}

In this section we extend a general contraction theorem for
recursive sequences $(Y_n)_{n\ge 0}$ of $d$-dimensional vectors as
developed in Neininger and R\"uschendorf \cite{neru4} to
families $(Y_t)_{t\ge 0}$ of $d$-dimensional vectors with continuous
parameter $t\in [0,\infty)$. For this, we assume that we have
\be \label{d-rec} Y_t \deq \sum_{r=1}^K A_r(t) Y^{(r)}_{T^{(t)}_r}
+b_t,\quad t\ge \tau_0, \ee
where $K$ is a positive integer and $\tau_0\ge 0$,
and we have $(Y^{(1)}_t)_t,\ldots,(Y^{(K)}_t)_t$,
$(A_1(t),\ldots,A_K(t),b_t,T^{(t)})_t$ independent,
where $T^{(t)}=(T^{(t)}_1,\ldots,T^{(t)}_K)$ is a
vector of random indices $T^{(t)}_r\in[0,t]$, the
$A_r(t)$ are random $d\times d$ matrices for
$r=1,\ldots,K$ and $b_t$ is a random $d$
dimensional vector. Finally, in (\ref{d-rec}), we
have that for each $t\ge 0$,  $Y_t$ and $Y^{(r)}_t$
are identically distributed for all $r=1,\ldots,K$.

We assume that all $Y_t$ as well as $A_r(t)$, $b_t$
and $T^{(t)}$ are defined on some probability space
$(\Omega,\cF,\mu)$, and that they are measurable
functions of $(t,\omega)$. (This is a technicality
to ensure that the sum in \eqref{d-rec} is
well-defined. Note, however, that the joint
distribution of $Y_t$ for different $t$ is
irrelevant.)

We introduce the normalized random vectors
\be\label{scald} X_t:=C_t^{-1/2}(Y_t-M_t),\quad t\ge 0, \ee
where $M_t\in\Rset^d$ and $C_t$ is a symmetric, positive definite
square matrix. We assume that $M_t$ and $C_t$ are measurable
functions of $t$; further restrictions on $M_t$ and $C_t$ will be
given in \convention{C}. The recurrence (\ref{d-rec}) implies a
recurrence for $X_t$,
\be\label{d-recsc} X_t \deq \sum_{r=1}^K A_r^{(t)}
X^{(r)}_{T^{(t)}_r} +b^{(t)},\quad t\ge \tau_0, \ee
 with  independence relations as in (\ref{d-rec}) and
\begin{align}
\label{modcoeff}
A_r^{(t)}=C_t^{-1/2}A_r(t)C_{T^{(t)}_r}^{1/2},
&&&
b^{(t)}=C_t^{-1/2}\left(b_t-M_t+\sum_{r=1}^K \bigpar{A_r(t)M_{T^{(t)}_r}}
\right).
\end{align}
As for the case with integer indexed vectors we establish a transfer
theorem of the following form: Appropriate convergence of the
coefficients $A_r^{(t)}\to A_r^\ast$, $b^{(t)}\to b^\ast$ implies
weak convergence of the quantities $X_t$ to a limit $X$. The
distribution ${\cL}(X)$ of $X$ is a fixed point of the limiting
equation obtained from (\ref{d-recsc}) by letting formally
$t\to\infty$:
 \be \label{limequ}X\deq  \sum_{r=1}^K A_r^\ast X^{(r)}
+b^\ast,\ee
where $(A_1^\ast,\ldots,A_K^\ast,b^\ast)$, $X^{(1)},\ldots,X^{(K)}$
are independent and $X^{(r)}\deq X$ for $r=1,\ldots, K$. To formalize
this
we introduce the map $T$ on the space $ \mathcal{M}^d$ of probability
measures on $\Rset^d$ by
\begin{align}\label{limmap0}
T: \mathcal{M}^d\to  \mathcal{M}^d,\quad \eta \mapsto \mathcal{L}
\left(\sum_{r=1}^K A_r^\ast Z^{(r)}+b^\ast\right),
\end{align}
where $(A_1^\ast,\ldots,A_K^\ast,b^\ast)$, $Z^{(1)},\ldots,Z^{(K)}$ are
independent and $\mathcal{L}(Z^{(r)})=\eta$ for $r=1,\ldots,K$. Then $X$
is a solution of  (\ref{limequ}) if and only if $ \mathcal{L}(X)$ is a
fixed point of $T$.

 We make use of Zolotarev's metric $\zeta_s$ with $0<s\le 3$. To ensure
 finiteness of the metric subsequently
 we make the following assumptions about the
scaling imposed in (\ref{scald}):

\begin{conventionC}
For $1<s\le 3$ we assume that $M_t=\E Y_t$.
For $2<s\le 3$ we assume that $\Cov(Y_t)$ is positive definite for all
$t\ge \tau_1$ with a $\tau_1\ge \tau_0$ and that $C_t=\Id_d$ for $0\le
t< \tau_1$ and
$C_t=\Cov(Y_t)$ for $t\ge \tau_1$.
\end{conventionC}

This convention implies that $X_t$ is centered for
$1<s\le 3$ and has $\Id_d$ as its covariance matrix
for $2< s\le 3$ and $t\ge \tau_1$. (For $0<s\le1$,
\convention C is void.)

\begin{thm}\label{theo1}
Let $0<s\le 3$ and let $(Y_t)_{t\ge 0}$ be a
process of random vectors satisfying \eqref{d-rec}
such that  $\|Y_t\|_s <\infty$ for every $t$.
Denote by $X_t$ the rescaled quantities in
\eqref{scald}, assuming \convention{C}. Assume that
$\| A_r^{(t)}\|_s<\infty$, $\|b^{(t)}\|_s<\infty$
and  $\sup_{0\le u\le t} \|X_u\|_s <\infty$ for
every $t\ge 0$, and
\begin{align}
&\left(A_1^{(t)},\ldots,A_K^{(t)},b^{(t)}\right)
\stackrel{\ell_s}{\longrightarrow}
\left(A_1^\ast,\ldots,A_K^\ast,b^\ast\right),
\label{con1bbb}\\
&\E \sum_{r=1}^K \|A_r^\ast\|_\op^s<1, \label{con2b}\\
&\E\left[\etta{\left\{T_r^{(t)}\le \tau\right\}}
\left\|A^{(t)}_r\right\|^s_\op\right] \to 0 \label{con3b}
\end{align}
for every $\tau>0$  and $r=1,\ldots,K$. Then $X_t$ converges in
distribution to  a limit $X$, and
\begin{equation}\label{zof}
\zeta_s(X_t,X)\to 0,\quad t\to\infty,
\end{equation}
where $\mathcal{L}(X)$ is the unique fixed point of $T$ given in
\eqref{limmap0} subject to $\|X\|_s<\infty$ and
\begin{equation}
  \label{cc1}
\begin{cases}
\E X =0 & \text{ for } 1<s\le 2,\\
\E X =0,\; \Cov(X)=\Id_d &\text{ for }   2<s\le 3.
\end{cases}
\end{equation}
\end{thm}
\begin{proof}
This proof is a continuous extension of the proof of  Theorem 4.1 in
Neininger and R\"uschendorf \cite{neru4} for the discrete time case.
The existence and uniqueness of the fixed point of $T$ subject to
(\ref{cc1}) is obtained as follows: For $1<s\le 3$ equation
(\ref{d-recsc}) implies $\E b^{(t)}=0$ for all $t>0$, thus by
(\ref{con1bbb}) we obtain $\E b^\ast =0$. For  $2<s\le 3$ equation
(\ref{d-recsc}) implies that for all $t\ge \tau_1$
\begin{align*}
\Id_d&=\Cov(X_t)\\
&=\E\left[b^{(t)}(b^{(t)})^\mathrm{tr}\right]+ \E \left[
\sum_{r=1}^K \Bigpar{\etta{\{T^{(t)}_r<\tau_1\}}A_r^{(t)}
\tilde{C}_{T^{(t)}_r} (A_r^{(t)})^\mathrm{tr}
 + \etta{\{T^{(t)}_r\ge\tau_1\}}A_r^{(t)} (A_r^{(t)})^\mathrm{tr}}\right],
\end{align*}
where $b^\mathrm{tr}$ denotes the transpose of a vector or matrix
and $\tilde{C}_t:=\Cov(X_t)$;
recall that $\tilde{C}_t=\Id$ when $t\ge\tau_1$.

By (\ref{con1bbb}), (\ref{con3b}) and
H\"older's inequality  this implies
\begin{align*}
\E\left[b^\ast(b^\ast)^\mathrm{tr}\right]+ \E \left[ \sum_{r=1}^K
A_r^\ast (A_r^\ast)^\mathrm{tr}\right]= \Id_d.
\end{align*}
Now, Corollary 3.4 in \cite{neru4} implies existence and uniqueness
of the fixed-point.\\
Since
\begin{equation}\label{512}
\E \sum_{r=1}^K \|A_r^{(t)}\|^s_\op \to \E \sum_{r=1}^K
\|A_r^\ast\|^s_\op = \xi <1
\end{equation}
there exist $\xi_+\in (\xi,1)$ and $\tau_2>\tau_1$ such that for all
$t\ge \tau_2$ we have
\begin{equation}\label{513}
\E \sum_{r=1}^K \|A_r^{(t)}\|^s_\op \le \xi_+ <1.
\end{equation}
 Now, we introduce the quantity
\begin{eqnarray}\label{qtdef}
Q_t:=\sum_{r=1}^K A_r^{(t)}
\left(\etta{\left\{T_r^{(t)}<\tau_2\right\}}X_{T^{(t)}_r}^{(r)}
 +\etta{\left\{T_r^{(t)}\ge\tau_2\right\}}X^{(r)}\right)
+b^{(t)},
\quad t\ge\tau_1,
\end{eqnarray}
where
$(A_1^{(t)},\ldots,A_K^{(t)},b^{(t)},T^{(t)}),X^{(1)},\ldots,X^{(K)},
(X^{(1)}_t),\ldots,(X^{(K)}_t)$ are independent with $X^{(r)}\sim X$
and $X^{(r)}_t \sim X_t$ for $r=1,\ldots,K$ and $t\ge 0$. Comparing
with (\ref{d-recsc}) we obtain that $Q_t$ is centered for $1<s\le 3$
and has the covariance matrix $\Id_d$ for $2<s\le 3$ and $t\ge
\tau_1$. Hence, $\zeta_s$ distances between $X_t$, $Q_t$ and $X$ are
finite for all $t\ge \tau_1$. The triangle inequality implies
\begin{eqnarray}\label{triang}
\Delta(t):=\zeta_s(X_t,X)\le \zeta_s(X_t,Q_t)+\zeta_s(Q_t,X).
\end{eqnarray}
As in the proof for the discrete case we obtain $\zeta_s(Q_t,X)\to
0$ as $t\to 0$, where we use that $\sup_{0\le
t\le\tau_2} \|X_t\|_s <\infty$.

The first summand of (\ref{triang}) requires a continuous analog
of the estimate in the discrete case.
Using the properties of the $\zeta_s$ metric, we obtain, for $t\ge \tau_1$,
\begin{eqnarray}
\zeta_s(X_t,Q_t)\le \E \sum_{r=1}^K \etta{\left\{T^{(t)}_r\ge
\tau_2\right\}}\left\|A^{(t)}_r\right\|^s_\op\Delta(T^{(t)}_r),
\end{eqnarray}
%{\tt !!! We need that  $t\mapsto \Delta(t)$ is measurable !!!}
%\marginpar{$\bullet\bullet$ GAP $\bullet\bullet$}\\
and,  with (\ref{triang}), and
$r_t:=\zeta_s(Q_t,X)$ it follows
\begin{eqnarray}\label{receq}
\Delta(t)\le \E \sum_{r=1}^K \etta{\left\{T^{(t)}_r\ge
\tau_2\right\}}\left\|A^{(t)}_r\right\|^s_\op\Delta(T^{(t)}_r) + r_t.
\end{eqnarray}
Now, we obtain $\Delta(t)\to 0$ in two steps, first
showing that $(\Delta(t))_{t\ge 0}$ is bounded and then, using the
bound, that $\Delta(t)\to 0$.

For the first step we introduce
\begin{equation}\label{joc}
\Delta^{\!\ast}(t):=\sup_{\tau_2\le u\le t} \Delta(u).
\end{equation}
We have $\Delta^{\!\ast}(t)<\infty$ for all $t\ge \tau_2$, since,
for $\tau_2\le u\le t$, we have $\zeta_s(X_u,X)\le C_s(\|X\|_s^s+
\|X_u\|_s^s) \le C_s(\|X\|_s^s+\sup_{\tau_2\le u\le
t}\|X_u\|_s^s)<\infty$ with a constant $C_s>0$, using \cite[Lemma
2]{Zolotarev76}. By definition, $\Delta^{\!\ast}$ is monotonically
increasing. With $R:=\sup_{t\ge \tau_2} r_t <\infty$ we obtain for
$\tau_2\le u\le t$,
from \eqref{receq}, \eqref{joc} and \eqref{513},
\begin{align*}
\Delta(u)&\le  \E \sum_{r=1}^K \etta{\left\{T^{(u)}_r\ge
\tau_2\right\}}\left\|A^{(u)}_r\right\|^s_\op\Delta^{\!\ast}(u)+ R\\
&\le \xi_+\Delta^{\!\ast}(t)+R.
\end{align*}
Hence, we obtain $\Delta^{\!\ast}(t)\le \xi_+\Delta^{\!\ast}(t)+R$, thus
 $\Delta^{\!\ast}(t)\le R/(1-\xi_+)$.
This implies
\begin{align}
\gDxoo:= \sup_{t\ge \tau_2}\Delta(t) \le \frac{R}{1-\xi_+}<\infty.
\end{align}

For the second step we denote
$L\=\limsup_{t\to\infty} \Delta(t)$.
%and assume, for a proof by contradiction,  that $L>0$.
For every $\eps>0$
there exists a $\tau_3>\tau_2$  such
that we have $\Delta(t)\le L+\varepsilon$ for all $t\ge \tau_3$.
Thus, from (\ref{receq}) we obtain
\begin{eqnarray*}
\Delta(t)\le \E \sum_{r=1}^K \etta{\left\{\tau_2\le T^{(t)}_r<
\tau_3\right\}}\left\|A^{(t)}_r\right\|^s_\op \gDxoo + \E   \sum_{r=1}^K
\etta{\left\{T^{(t)}_r\ge
\tau_3\right\}}\left\|A^{(t)}_r\right\|^s_\op    (L+\varepsilon)
 + r_t
\end{eqnarray*}
and letting $t\to\infty$ we obtain by (\ref{con3b}) and \eqref{512}
\begin{eqnarray*}
L\le \xi (L+\varepsilon).
\end{eqnarray*}
If $L>0$, this is a contradiction for
$0<\varepsilon<L(1-\xi)/\xi$.
Hence, we have $L=0$. This proves \eqref{zof}.
Finally, recall that convergence in $\zeta_s$ implies weak
convergence.
\end{proof}

As a corollary we formulate a univariate central limit theorem
that corresponds to Neininger and R\"uschendorf
\cite[Corollary 5.2]{neru4}
for the discrete time case. For this we assume that there are
expansions, as $t\to\infty$,
\begin{align}\label{expan}
\E Y_t=f(t)+o(g^{1/2}(t)),&&& \Var(Y_t)=g(t)+o(g(t))
\end{align}
with functions $f:[0,\infty)\to \Rset$,
$g:[0,\infty)\to [0,\infty)$, with
\begin{align}\label{regul}
 \sup_{u\le t} |f(u)|<\infty \text{ for every $t>0$},
&&&
\lim_{t\to\infty}g(t)=\infty,
&& \sup_{u\le t} g(u)= O(g(t)).
\end{align}
Thus, for some constant $C\ge1$, $g(u)\le Cg(t)$ when $0\le u\le t$.

Then the following central limit law holds:
\begin{cor} \label{coro}
Let $2<s\le3$ and let $Y_t$, ${t\ge0}$, be given
$s$-integrable, univariate random  variables satisfying \eqref{d-rec}
  with $A_r(t)=1$ for all $r=1,\ldots,K$ and
  $t\ge 0$.
Assume that $\sup_{u\le t} \E|Y_u|^s<\infty$ for
every $t$,  and that the mean and variance of\/
  $Y_t$ satisfy
\eqref{expan} with  \eqref{regul}.
If, as \ttoo,
\begin{align}\label{con1}
\left(\sqrt{\frac{g(T^{(t)}_1)}{g(t)}},
\ldots,\sqrt{\frac{g(T^{(t)}_K)}{g(t)}}\right)
&\stackrel{\ell_s}{\longrightarrow} (A_1^\ast,\ldots,A_K^\ast),\\
\frac{1}{g^{1/2}(t)}\left(b_t-f(t)+\sum_{r=1}^K
  f(T^{(t)}_r)\right)
&\stackrel{\ell_s}{\longrightarrow} 0,\label{con1b}
\end{align}
and furthermore
\begin{align}\label{con2}
\sum_{r=1}^K (A_r^\ast)^2=1,
&&&
\Prob\left(\bigcup_{r=1}^K \{A_r^\ast=1\}\right)<1,
\end{align}
then
\begin{align}\label{con3}
\frac{Y_t-f(t)}{g^{1/2}(t)}\dto {\cN}(0,1).
\end{align}
\end{cor}

\begin{proof}
We begin by replacing $g(t)$ by $\max(g(t),1)$;
by \eqref{regul}, this does not affect $g(t)$ for large $t$, and it is
easy to see that
\eqref{expan}, \eqref{regul}, \eqref{con1}, \eqref{con1b} still hold.
We may thus assume that $g(t)\ge1$ for every $t$.

Denote $M_t\=\E Y_t$ and $\gsst\=\Var(Y_t)$. By \eqref{expan},
$\gsst/g(t)\to1$.
(All unspecified limits are as \ttoo.)
Choose $\tau_1\ge\tau_0$ such that $\tfrac14 g(t)\le\gsst\le4g(t)$
for $t\ge\tau_1$.
Let, as in \convention{C}, $C_t\=1$ for $t<\tau_1$ and
$C_t\=\gsst$ for $t\ge\tau_1$, and write $\tgs_t\= C_t\qh$ and
$\eps(t)\=\tgs_t/g(t)\qh-1=\bigpar{C_t/g(t)}\qh-1$.
For $t\ge\tau_1$, $\eps(t)=\bigpar{\Var Y_t/g(t)}\qh-1$, so by
\eqref{expan},
\begin{equation}
  \label{s1}
\eps(t)\to0
\qquad\text{as \ttoo}.
\end{equation}
Further, $C_t/g(t)=1/g(t)\le1$ for $t<\tau_1$, while $C_t/g(t)=\gsst/g(t)\le 4$
for $t\ge\tau_1$. Hence $|\eps(t)|\le1$ for all $t$.
With \eqref{modcoeff} and $A_r(t)=1$ we have, for $t\ge\tau_1$,
\begin{align}
  \art
=\frac{\tgs_{\trt}}{\gst}
=\frac{(1+\eps(\trt))g(\trt)\qh}{\gst},
\label{slask}
\\
\bt=\gst^{-1}\left(b_t-M_t+\sum_{r=1}^K M_{\trt}
\right).
\label{bt}
\end{align}
Since $g(\trt)\le Cg(t)$ by \eqref{regul}, we have, for $t\ge\tau_1$,
\begin{equation}\label{s2}
\Norm{  \art-\frac{g(\trt)\qh}{\gst}}_s
=
\Norm{  \eps(\trt)\frac{g(\trt)\qh}{\gst}}_s
\le
\sup_{u \le t}\lrabs{\eps(u)\frac{g(u)\qh}{\gst}}
.
%\le C \Norm{  \eps(\trt)\frac{g(t)\qh}{\gst}}_s
%\le 2C \Norm{  \eps(\trt)}_s.
\end{equation}
For any $\gd>0$, there exists, by \eqref{s1},
$\taud\ge\tau_1$ such that  $|\eps(t)|\le\gd$ when
$t\ge\taud$. Thus, if $\taud\le u\le t$, then
\begin{equation*}
\lrabs{\eps(u)\frac{g(u)\qh}{\gst}}
\le
\gd\,\frac{C g(t)\qh}{\gst}
\le 2C\gd.
\end{equation*}
On the other hand, if $u\le\taud$, then
\begin{equation*}
\lrabs{\eps(u)\frac{g(u)\qh}{\gst}}
\le
\frac{C g(\taud)\qh}{\gst}
\to0
\end{equation*}
as \ttoo. Hence,
$\sup_{u\le t}|\eps(u)\xfrac{g(u)\qh}{\gst}| \le 2C\gd$ for
sufficiently large $t$. Since $\gd>0$ is arbitrary, it follows that
the right hand side of \eqref{s2} tends to 0 as \ttoo, and thus
\eqref{s2} yields
\begin{equation}\label{s2+}
\Norm{  \art-\frac{g(\trt)\qh}{\gst}}_s
\to0.
\end{equation}
Since $g(t)\qh/\gst\to1$, \eqref{con1} yields
$g(\trt)\qh/\gst\lsto \axr$, which combined with \eqref{s2+} yields
$\art\lsto\axr$, jointly for $r=1,\dots,k$.

Next, for any $\eps>0$, there exists by \eqref{expan} $\taue\ge\tau_1$ such
that
$|M_t-f(t)|\le \eps g(t)\qh$ if $t\ge\taue$. Consequently,
if $\trt\ge\taue$, then
\begin{equation*}
  |M_{\trt}-f(\trt)|
\le \eps g(\trt)\qh
\le C\eps g(t)\qh.
\end{equation*}
Since $\sup_{u\le\taue}|M_u|$ and $\sup_{u\le\taue}|f(u)|$ are finite,
the same estimate holds for $\trt<\taue$ too, provided $t$ is large.
Consequently, $|M_{\trt}-f(\trt)|/g(t)\qh\le C\eps$ if $t$ is large
enough. It follows that
$\norm{M_{\trt}-f(\trt)}_s/g(t)\qh\to0$ as \ttoo,
so by \eqref{bt}, \eqref{con1b} and \eqref{expan},
$\bt\lsto0$.

We apply Theorem \ref{theo1} with $2<s\le 3$; we
have shown that \eqref{con1bbb} holds with
$b^\ast=0$. The two assumptions in (\ref{con2}) and
$s>2$ ensure that we have $\E \sum_{r=1}^K
|A_r^\ast|^s<1$. Finally,  by \eqref{s2+}, for
every $\tau$ and $r$,
\begin{equation*}
  \begin{split}
\Norm{\etta{\left\{T_r^{(t)}\le \tau\right\}}\art}_s
&\le
\Norm{\etta{\left\{T_r^{(t)}\le \tau\right\}}\frac{g(\trt)\qh}{\gst}}_s
+
\Norm{\art-\frac{g(\trt)\qh}{\gst}}_s
\\&
\le
\frac{C g(\tau)\qh}{\gst}+ o(1)
\to0.
  \end{split}
\end{equation*}

Now, Theorem \ref{theo1} implies $(Y_t-M_t)/\sigma_t\dto X$, where
${\cL}(X)$ is characterized by $\|X\|_s<\infty$, $\E X=0$,
$\Var(X)=1$, and
\begin{align}
X\deq \sum_{r=1}^K A_r^\ast X^{(r)},
\end{align}
with assumptions as in (\ref{limequ}). Since $\sum_{r=1}^K (A_r^\ast)^2=1$
this is solved by
${\cL}(X)={\cN}(0,1)$.
Consequently,
\begin{align*}
\frac{Y_t-M_t}{\sigma_t}\dto {\cN}(0,1),
\end{align*}
which, in view of (\ref{expan}), implies the assertion.
\end{proof}

The following theorem covers cases where the previous central limit
theorem of Corollary \ref{coro} fails due to the appearance of
periodic behavior. For this we assume that there is an expansion of
the mean, as $t\to\infty$,
\begin{align}\label{expan2}
\E Y_t=f(t)+\Re(\gamma t^\lambda)+o(t^\sigma),
\end{align}
with a function $f:[0,\infty)\to \Rset$, $\gamma \in
\mathbb{C}\setminus \{0\}$, and $\lambda \in \mathbb{C}$ with
$\sigma:=\Re(\lambda)>0$. We denote
\begin{align}
 &A_r^{(t)}\=\left(\frac{T_r^{(t)}}{t}\right)^\lambda,\quad
  r=1,\ldots,K,\label{06r2}\\
  &b^{(t)}\=\frac{1}{t^\sigma}\left(b_t-f(t)+\sum_{r=1}^K
  f(T^{(t)}_r)\right).\label{06r1}
\end{align}
Note that $\art$ in general is complex, while
$b^{(t)}$ is real.

% Furthermore we assume that
%\begin{equation}
%  \label{regul2b}
%T^{(t)}_r\pto\infty \qquad\text{as \ttoo,\quad $r=1,\dots,K$}.
%\end{equation}

\begin{thm} \label{gthper}
Let $Y_t$, ${t\ge0}$, be given square-integrable, univariate random
variables satisfying \eqref{d-rec}
  with $A_r(t)=1$ for all $r=1,\ldots,K$ and
  $t\ge 0$.
Assume that $\sup_{u\le t} \E|Y_u|^2<\infty$ for every $t>0$ and
that the mean of
  $Y_t$ satisfies
\eqref{expan2} with $\gl=\gs+\ii\tau$ and $\gs>0$,
 and some locally bounded function $f(t)$. If, as
\ttoo,
\begin{align}
(A_1^{(t)},\ldots,A_K^{(t)}) \stackrel{\ell_2}{\longrightarrow}
(A_1^\ast,\ldots,A_K^\ast)
&&&\text{and}&&
\|b^{(t)}\|_2\to 0,\label{con1bneu}
\end{align}
and furthermore
\begin{align}\label{con2neu}
\E \sum_{r=1}^K |A_r^\ast|^{2}<1,
\end{align}
then, as \ttoo,
\begin{align}\label{con3neu}
\ell_2\left(\frac{Y_t-f(t)}{t^\sigma},\Re\left(X e^{\ii\tau\ln
t}\right)\right) \to 0,
\end{align}
where $\mathcal{L}(X)$ is the unique fixed point in
$\mathcal{M}^\mathbb{C}_2(\gamma)$ of
\begin{align}\label{limmap0neu}
T: \mathcal{M}^\mathbb{C}\to  \mathcal{M}^\mathbb{C},\quad \eta
\mapsto \mathcal{L} \left(\sum_{r=1}^K A_r^\ast Z^{(r)}\right),
\end{align}
where $(A_1^\ast,\ldots,A_K^\ast)$, $Z^{(1)},\ldots,Z^{(K)}$ are
independent and $\mathcal{L}(Z^{(r)})=\eta$ for $r=1,\ldots,K$.
\end{thm}
\begin{proof}
We extend an approach based on the contraction
method from Fill and Kapur \cite{fika04}. We may
assume that $\tau_0\ge1$.

First,  for technical convenience we show that we
further may assume $Y_t=0$ and $f(t)=0$ for $0\le
t\le 1$.
%(An alternative would be to replace $t$ by $t\vee1$ and $T_r^{(t)}$
Let $(Y^{(r)*}_{t})_t$, $r=1,\dots,K$, be another set of copies of $(Y_t)_t$, independent of each other and of everything else.
We may replace $Y^{(r)}_{t}$ in \eqref{d-rec} by
$Y^{(r)}_{t}\ettax{t\ge1}+Y^{(r)*}_{t}\ettax{t<1}$, which has the same distribution and independence properties.
Hence $Y^{(r)}_{t}\ettax{t\ge1}$ satisfies \eqref{d-rec} (for $t\ge\tau_0\ge1$)
with $b_t$ replaced by
$\tilde b_t:=b_t+\sum_r Y^{(r)*}_{T_r^{(t)}}\ettax{T_r^{(t)}<1}$.
This replaces $b^{(t)}$ by $\tilde b^{(t)}$ with
\begin{equation*}
\lrabs{\tilde b^{(t)}-b^{(t)}}
\le
t^{-\gs}\sum_r\lrabs{ Y^{(r)*}_{T_r^{(t)}}}\ettax{T_r^{(t)}<1}
\end{equation*}
so
$\norm{\tilde b^{(t)}-b^{(t)}}_2=O(t^{-\gs})$ and \eqref{con1bneu} still holds.
We may thus consider $Y^{(r)}_{t}\ettax{t\ge1}$ instead, and thus we may assume that
$Y^{(r)}_{t}=0$ when $t<1$.
Similarly, we may assume that $f(t)=0$ for $t<1$, changing
$b^{(t)}$ by $O(t^{-\gs})$.

With $X_t\=(Y_t-f(t))/t^\sigma$ for $t>0$ and
$X_0:=0$ we obtain
\begin{align}\label{06r3}
X_t\eqd \sum_{r=1}^K \left(\frac{T_r^{(t)}}{t}\right)^\sigma
X_{T_r^{(t)}}^{(r)}+b^{(t)},
\qquad t\ge\tau_0,
\end{align}
with $b^{(t)}$ as given in (\ref{06r1}).

Next we prove that the restriction of $T$ defined in
(\ref{limmap0neu}) to $\mathcal{M}^\mathbb{C}_2(\gamma)$  maps into
$\mathcal{M}^\mathbb{C}_2(\gamma)$ and is Lipschitz in $\ell_2$ with
Lipschitz constant bounded by $\bigpar{\E \sum_{r=1}^K
|A_r^\ast|^{2}}\qh<1$.
%; see R\"osler and R\"uschendorf \cite[Lemma 1]{roerue}.

Note that (\ref{con2neu}) implies $\|A_r^\ast\|_2<\infty$  for all
$r=1,\ldots,K$. This implies that $T(\eta)$ has a finite second
moment for all $\eta\in \mathcal{M}^\mathbb{C}_2$. Next we claim
that $\sum_{r=1}^K \E A_r^\ast =1$. This implies that $T(\eta)$ has
mean $\gamma$ for all $\eta\in\mathcal{M}^\mathbb{C}_2(\gamma)$. To
prove $\sum_{r=1}^K \E A_r^\ast =1$,  note that (\ref{expan2})
implies $\E X_t=\Re(\gamma t^{\ii\tau})+o(1)$ as $t\to\infty$. On the
other hand, the right hand side of (\ref{06r3}) has mean, using $\E
b^{(t)}\to 0$,
\begin{align*}
\sum_{r=1}^K \E \left[ \left(\frac{T_r^{(t)}}{t}\right)^\sigma
\Re\bigpar{\gamma(T_r^{(t)})^{\ii\tau}}  \right]+o(1) &= \Re\left(
\gamma\sum_{r=1}^K \E
\frac{(T^{(t)}_r)^\lambda}{t^\sigma}\right)+o(1)\\
&= \Re\left( \gamma t^{\ii\tau}\sum_{r=1}^K \E
\left(\frac{T^{(t)}_r}{t}\right)^\lambda\right)+o(1)\\
&= \Re\left( \gamma t^{\ii\tau}\sum_{r=1}^K \E A_r^\ast\right)+o(1),
\end{align*}
where we also used that $\E (T^{(t)}_r/t)^\lambda \to \E A_r^\ast$,
see (\ref{con1bneu}). Hence, together we obtain, as $t\to\infty$,
\begin{align}
\Re(\gamma t^{\ii\tau})+o(1) = \Re\left( \gamma t^{\ii\tau}\sum_{r=1}^K
\E A_r^\ast\right)+o(1).
\end{align}
Thus, $\gamma\neq 0$ yields $\sum_{r=1}^K \E
A_r^\ast =1$. For the bound on the Lipschitz
constant in $\ell_2$ of $T$ restricted to
$\mathcal{M}^\mathbb{C}_2$ see R\"osler and
R\"uschendorf \cite[Lemma 1]{roerue} and Fill and
Kapur \cite{fika04}: For $\mu,\nu\in
\mathcal{M}^\mathbb{C}_2$ choose
$(Z^{(1)},W^{(1)}),\ldots,(Z^{(K)},W^{(K)})$ as
identically distributed vectors of optimal
couplings of $\mu$ and $\nu$ and such that
$(Z^{(1)},W^{(1)}),\ldots,(Z^{(K)},W^{(K)})$,
$(A_1^\ast,\ldots,A_K^\ast)$ are independent. Then
we have
\begin{align*}
\ell_2^2(T(\mu),T(\nu)) &=
\ell_2^2\left(\sum_{r=1}^K A_r^\ast Z^{(r)},
\sum_{r=1}^K A_r^\ast W^{(r)}\right)\\
&\le \E \left|\sum_{r=1}^K A_r^\ast(Z^{(r)}-W^{(r)})\right|^2\\
&=\E \left(\sum_{r=1}^K
|A_r^\ast|^2|Z^{(r)}-W^{(r)}|^2 +\sum_{r\neq s}
A_r^\ast(Z^{(r)}-W^{(r)})\overline{A_s^\ast
(Z^{(s)}-W^{(s)})} \right)
\\
&=\E\sum_{r=1}^K  |A_r^\ast|^2 \ell_2^2(\mu,\nu)+0
\\
&=\sum_{r=1}^K \E |A_r^\ast|^2 \ell_2^2(\mu,\nu).
\end{align*}
Altogether we obtain that $T$ has a unique fixed point
$\mathcal{L}(X)$ in $\mathcal{M}^\mathbb{C}_2(\gamma)$.

The fixed point property of $\mathcal{L}(X)$
implies
\begin{eqnarray}\label{06nn1}
\frac{1}{t^\sigma}\Re\left(t^\lambda X\right)
\eqd\frac{1}{t^\sigma}\Re\left(\sum_{r=1}^K t^\lambda A_r^\ast
X^{(r)}\right).
\end{eqnarray}
where $(A_1^\ast,\ldots,A_K^\ast)$,
$X^{(1)},\ldots,X^{(b)}$ are independent and
$\mathcal{L}(X^{(r)})=\mathcal{L}(X)$ for
$r=1,\ldots,K$. We may assume, e.g.\ by taking
optimal couplings, that
$\norm{A_r^{(t)}-A_r^*}_2\to0$ as \ttoo. We choose
$X_t^{(r)}$ as optimal couplings to
%$t^{-\gs}\Re(t^\gl X^{(r)})$
$\Re(t^{\ii\tau} X^{(r)})$  (with the right
distribution, \ie{} the distribution of $X_t$) for
$t\ge 0$ and $r=1,\ldots,K$.  Clearly, we may
assume that, as required, $X_t^{(r)}$,
$r=1,\dots,K$, are independent of each other and of
$(T^{(t)},b_t)_t$.

 We denote, for $t>0$,
\begin{eqnarray*}
\Delta(t)\=\ell_2\left(\frac{Y_t-f(t)}{t^\sigma},\Re\left(X
e^{\ii\tau\ln
t}\right)\right)=\ell_2\left(X_t,\frac{1}{t^\sigma}\Re\left(t^\lambda
X\right)\right).
\end{eqnarray*}
Using (\ref{06r3}) and (\ref{06nn1}) we obtain, for
$t\ge\tau_0$,
\begin{align}
\Delta(t)&=\ell_2\left(\sum_{r=1}^K
\left(\frac{T_r^{(t)}}{t}\right)^\sigma
X_{T_r^{(t)}}^{(r)}+b^{(t)},\frac{1}{t^\sigma}\Re\left(\sum_{r=1}^K
t^\lambda A_r^\ast X^{(r)}\right)\right)\nonumber\\
&\le \left\| \sum_{r=1}^K \left(
\left(\frac{T_r^{(t)}}{t}\right)^\sigma
X_{T_r^{(t)}}^{(r)}-\frac{1}{t^\sigma}\Re\left(t^\lambda A_r^\ast
X^{(r)}\right)\right)\right\|_2+\left\|b^{(t)}\right\|_2\nonumber\\
&\le \left\| \sum_{r=1}^K \left(
\left(\frac{T_r^{(t)}}{t}\right)^\sigma
X_{T_r^{(t)}}^{(r)}-\frac{1}{t^\sigma}\Re\left((T_r^{(t)})^\lambda
X^{(r)}\right)\right)\right\|_2+\left\|b^{(t)}\right\|_2\nonumber\\
&\qquad{}+ \left\| \sum_{r=1}^K
\left(\frac{1}{t^\sigma}\Re\left((T_r^{(t)})^\lambda
X^{(r)}\right)-\frac{1}{t^\sigma}\Re\left(t^\lambda A_r^\ast
X^{(r)}\right)\right)\right\|_2\label{06r5}
.
\end{align}
By (\ref{con1bneu}) and (\ref{06r2}) the second and
third of the three latter summands tend to zero as
\ttoo. We abbreviate
\begin{align}
W_r^{(t)}\=\left(\frac{T_r^{(t)}}{t}\right)^\sigma
X_{T_r^{(t)}}^{(r)}-\frac{1}{t^\sigma}\Re\left((T_r^{(t)})^\lambda
X^{(r)}\right).
\end{align}
Hence, (\ref{06r5}) implies
\begin{align}
\Delta(t)&\le\left(\E \left(\sum_{r=1}^K
W_r^{(t)}\right)^2\right)^{1/2}+o(1)\nonumber\\
&= \left(\E \sum_{r=1}^K (W_r^{(t)})^2
+ \E \sum_{\substack{r,s=1\\ r\neq s}}^K
  W_r^{(t)} W_s^{(t)}\right)^{1/2}+o(1). \label{06r9}
\end{align}
By the definition of $\Delta(t)$ and the fact that
$\bigpar{X^{(r)}_t,\Re(t^{\ii\tau} X^{(r)})}$  are
optimal couplings for all $t>0$ and $r=1,\ldots,K$
we obtain
\begin{align}\label{06r8}
\E (W_r^{(t)})^2 =\E\left[
\left(\frac{T_r^{(t)}}{t}\right)^{2\sigma}\Delta^2(T_r^{(t)})\right].
\end{align}
{}From (\ref{expan2}) we obtain
\begin{align*}
\E X_t = \frac{1}{t^\sigma}\Re(\gamma t^\lambda)+R(t),
\qquad t>0,
\end{align*}
with $R(t)\to 0$ as \ttoo. Since $\E
X^{(r)}=\gamma$ and by the independence conditions
we obtain $\E W_r^{(t)}=\E[(T^{(t)}_r/t)^\sigma
R(T_r^{(t)})]$ and, for $r\neq s$,
\begin{align*}
\E[W_r^{(t)}W_s^{(t)}]&=\E\left[\left(\frac{T_r^{(t)}}{t}
\frac{T_s^{(t)}}{t}\right)^\sigma R(T_r^{(t)})R(T_s^{(t)})\right].
\end{align*}
Splitting the latter integral into the events
$\{T^{(t)}_r\le \taux_1 \text{ or } T^{(t)}_s\le
\taux_1\}$ and $\{T^{(t)}_r> \taux_1 \text{ and }
T^{(t)}_s>\taux_1\}$ for some $\taux_1>0$ we
obtain, for every $\taux_1>0$,
\begin{align*}
\lrabs{\E[W_r^{(t)}W_s^{(t)}]}\le \left(\frac{\taux_1}{t}\right)^\sigma \norm{R}_\infty^2
+\sup_{u\ge \taux_1}R^2(u),
\end{align*}
where $\|R\|_\infty:=\sup_t|R(t)|<\infty$.  {}From
this we obtain first, letting $\ttoo$,
$\limsup_{\ttoo}\lrabs{\E[W_r^{(t)}W_s^{(t)}]}\le
\sup_{u\ge \taux_1}R^2(u)$, and then, letting
$\taux_1\to\infty$,
\begin{align}\label{06r6}
\E[W_r^{(t)}W_s^{(t)}]\to 0
\qquad \text{as }t\to\infty.
\end{align}

Now, (\ref{06r9}), (\ref{06r8}), and (\ref{06r6}) imply, for $t>\tau_0$,
\begin{eqnarray}
 \Delta(t)\le \left( \E \left[\sum_{r=1}^K
 \left(\frac{T_r^{(t)}}{t}\right)^{2\sigma}
 \Delta^2(T_r^{(t)})\right]+R_1(t)\right)^{1/2}+R_2(t),\label{06n10}
\end{eqnarray}
with $R_1(t),R_2(t) \to 0$ as \ttoo.
%Although we assumed $t\ge\tau_0$ above, \eqref{06n10} holds for all $t>0$ if we choose
%$R_2(t)$ large enough for $t\le\tau_0$. Moreover,
%We may thus assume that $R_1(t)$ and $R_2(t)$ in \eqref{06n10} are bounded,
%since we otherwise may replace them by suitable constants on some interval $[0,\tau_0']$.

We first show that $\|\Delta\|_\infty<\infty$.
Define $\Delta^\ast (t) \=\sup_{0< u\le t}
\Delta(u)$. By the assumptions $\sup_{0\le u\le
t}\E |Y_u|^2<\infty$ and $\sup_{0\le u\le t}
|f(u)|<\infty$, together with $Y_u=0$ and $f(u)=0$
for $u\le1$, we have $\Delta^\ast(t)<\infty$ for
all $t>0$. Let $t_1\ge\tau_0$ be such that
$|R_1(t)|<1$ and $|R_2(t)|<1$ for $t\ge t_1$. Then
with (\ref{06n10}) we obtain, for $t\ge t_1$,
\begin{equation*}
 \Delta(t)\le \left( \E \left[\sum_{r=1}^K
 \left(\frac{T_r^{(t)}}{t}\right)^{2\sigma}
 (\Delta^\ast)^2(t)\right]+1\right)^{1/2}+1.
\end{equation*}
By \eqref{06r2}, \eqref{con1bneu} and (\ref{con2neu}) there exists a $t_2\ge t_1$
such that for all $t\ge t_2$ we have $\E \sum_{r=1}^K
 (T_r^{(t)}/t)^{2\sigma}\le \xi<1$. Thus, for all
 $t\ge t_2$ we obtain,
 with $\sqrt{a+b}\le \sqrt{a}+\sqrt{b}$ for $a,b\ge 0$,
\begin{equation*}
 \Delta(t)\le \sqrt{\xi}\gD^\ast(t)+ 2,
\end{equation*}
and thus
\begin{equation*}
 \Delta^\ast(t)\le \sqrt{\xi}\gD^\ast(t)+ 2 +\gD^\ast(t_2),
\end{equation*}
which implies
$\|\Delta\|_\infty
\le\bigpar{2+\gD^\ast(t_2)}/(1-\sqrt\xi)
<\infty$.

In a second step we show that $\Delta(t)\to 0$ as
\ttoo. For this we assume that $L:=\limsup_{t \to
\infty} \Delta(t)>0$. Let $\varepsilon>0$.  There
exists a $t_3\ge t_2$ such that for all $t\ge t_3$
we have $\Delta(t)\le L+\varepsilon$. Then
(\ref{06n10}) implies
\begin{eqnarray*}
\lefteqn{ \Delta(t)}\\
&\le& \left( \E \left[\sum_{r=1}^K
 \left(\frac{T_r^{(t)}}{t}\right)^{2\sigma}
 \left(\etta{\{T^{(t)}_r<t_3\}}+\etta{\{T^{(t)}_r\ge  t_3\}}\right)
 \Delta^2(T_r^{(t)})\right]+R_1(t)\right)^{1/2}+R_2(t)\\
 &\le&
 \left(\sum_{r=1}^K\left(\frac{t_3}{t}\right)^{2\sigma}\|\Delta\|_\infty^2
 +\xi(L+\varepsilon)^2+R_1(t)\right)^{1/2}+R_2(t).
\end{eqnarray*}
Hence, $t\to \infty$ implies
\begin{eqnarray*}
 L\le \sqrt{\xi}(L+\varepsilon),
\end{eqnarray*}
which if $L>0$ is a contradiction if we choose
%$\varepsilon<(1-\sqrt{\xi})L/\sqrt{\xi}$.
$\varepsilon$ small enough.
Consequently, we
have $L=0$ yielding the assertion.
\end{proof}

\section{Proof of \refT{T1}} \label{Sproof}
 In this section we prove \refT{T1}. The statements on
 mean and variance of $N(x)$ are proved in Section \ref{Smean}.
 It remains to identify the asymptotic distribution of $N(x)$
 Note that recurrence (\ref{a5}) for $N(x)$ is covered by the
 general recurrence for $Y_t$ in (\ref{d-rec}) by making the
 choices $d=1$, $K=b$, $\tau_0=1$, $A_r(t)=1$, $T_r^{(t)}=V_rt$
 and $b_t=1$ for all $r=1,\ldots,K$
 and $t\ge \tau_0$.

 We consider the three cases (i) -- (iii) appearing in
 \refT{T1} separately:\\

\noindent Case (i):
We have $\E N(x)=\alpha^{-1}x+o(\sqrt{x})$ and
 $\Var(N(x))\sim \beta x$ with $\beta>0$. We apply
 \refC{coro} with the choices $f(t)=\alpha^{-1} t$ and
 $g(t)=\beta t$. The conditions (\ref{expan}) and (\ref{regul})
 are satisfied.
%Since $V_r$ is absolutely continuous on $[0,1)$ we obtain $T^{(t)}_r= V_rt\pto \infty$ as $t\to\infty$.
We have
 $\sup_{u\le t} \E |Y_u|^s<\infty$ for $s=3$ by Lemma \ref{Lmom}.
 Condition (\ref{con1}) is satisfied with $A_r^\ast=\sqrt{V_r}$ for
 $r=1,\ldots,K$, condition (\ref{con1b}) is trivially satisfied, and
 we have (\ref{con2}). Hence, Corollary \ref{coro} applies and
 yields
\begin{eqnarray*}
 \frac{N(x)-\alpha^{-1}x}{\sqrt{\beta x}}\dto \mathcal{N}(0,1),
\end{eqnarray*}
which is the assertion.\\

\noindent Case (ii):
We have $\E N(x)=\alpha^{-1}x+O(\sqrt{x})$ and
 $\Var(N(x))\sim \beta x\ln x$ with $\beta>0$. We apply
\refC{coro} with the choices $f(t)=\alpha^{-1} t$ and
 $g(t)=\beta t\ln t$. Now we have
 $g(T^{(t)}_r)/g(t)=V_r+V_r\ln(V_r)/\ln t$, hence we obtain, since
 $x\mapsto x\ln x$ is bounded on $[0,1]$,
\begin{eqnarray*}
\left(\sqrt{\frac{g(T^{(t)}_1)}{g(t)}},
\ldots,\sqrt{\frac{g(T^{(t)}_K)}{g(t)}}\right)
&\stackrel{\ell_3}{\longrightarrow} (A_1^\ast,\ldots,A_K^\ast),
\end{eqnarray*}
with $A_r^\ast=\sqrt{V_r}$ for
 $r=1,\ldots,K$. All conditions of Corollary \ref{coro} are
 satisfied as in case (i) and we obtain
\begin{eqnarray*}
 \frac{N(x)-\alpha^{-1}x}{\sqrt{\beta x\ln x}}\dto
 \mathcal{N}(0,1).
\end{eqnarray*}

\noindent Case (iii):
\refT{gthper} can be applied and yields
\begin{align*}
\ell_2\left(\frac{\nx  -\gai x}{x^{\Re\gl_2}}, \Re\bigpar{\Xi
e^{\ii\Im\gl_2\ln x}}\right) \to 0
\end{align*}
as $x\to \infty$.

Here, we give a simplified version of the proof of \refT{gthper} for
the special recurrence (\ref{a5}) which yields also the rate of
convergence stated in \refT{T1}.

The restriction of $T$ defined in (\ref{limmap}) to
$\mathcal{M}^\mathbb{C}_2(\nu)$ is Lipschitz in $\ell_2$ with
Lipschitz constant bounded by  $\bigpar{\E \sum_{r=1}^b
V_r^{2\Re(\lambda_2)}}\qh$; cf.~the first part of the proof of Theorem
\ref{gthper}. Recall $\lambda_2=\sigma+\ii\tau$ with real
$\sigma,\tau$. Hence, $\sigma>1/2$  implies that $T$ has a unique
fixed point $\mathcal{L}(\Xi)$ in $\mathcal{M}^\mathbb{C}_2(\nu)$.

 For $X_t\=N(t)-\alpha^{-1} t$ we obtain with
(\ref{a5})
\begin{align}
X_t\deq \sum_{r=1}^b  X_{V_r t}^{(r)} + 1,
\end{align}
where  $X_{t}^{(r)}$ are independent distributional copies of $X_t$
also independent of $(V_1,\ldots,V_b)$. With the fixed point
property of $\Xi$  we have
\begin{align*}
 t^{\sigma}\Re \left(\Xi e^{\ii\tau \ln x}\right)=\Re(t^{\lambda_2}\Xi)\deq \Re\left(\sum_{r=1}^b
 (V_rt)^{\lambda_2} \Xi^{(r)}\right),
\end{align*}
where $(V_1,\ldots,V_b)$,
$\Xi^{(1)},\ldots,\Xi^{(b)}$ are independent and
$\mathcal{L}(\Xi^{(r)})=\mathcal{L}(\Xi)$ for
$r=1,\ldots,b$. We choose $X_t^{(r)}$ as optimal
couplings to $\Re(t^{\gl_2}\Xi^{(r)})$  for $t\ge
0$ and $r=1,\ldots,b$ and denote
$\Delta(t):=\ell_2(X_t,\Re(t^{\lambda_2}\Xi))$.
Note that in the definition of $X_t$ we did not
rescale by $t^\sigma$, hence we have to show
$\Delta(t)=O(t^\kappa)$.

With $W^{(t)}_r:=X_{V_rt}^{(r)}-\Re((V_rt)^{\lambda_2}\Xi^{(r)})$ we
obtain, for $t\ge1$,
\begin{align*}
 \Delta(t)&=\ell_2\left(\sum_{r=1}^b X_{V_rt}^{(r)}+1,
\sum_{r=1}^b \Re\left((V_rt)^{\lambda_2}\Xi^{(r)}\right)\right)\\
 &\le \Bigg\{\E \Bigg(\sum_{r=1}^b W^{(t)}_r\Bigg)^2\Bigg\}^{1/2}+1\\
 &=\Bigg\{\sum_{r=1}^b \E (W^{(t)}_r) ^2+
 \sum_{\substack{r,s=1 \\ r\neq  s}}^b\E[ W^{(t)}_rW^{(t)}_s]   \Bigg\}^{1/2}+1.
\end{align*}
Conditioning on $(V_1,\ldots,V_b)$ yields $\E
(W^{(t)}_r)^2=\E\Delta^2(V_rt)$.
 From $\E N(t) =\alpha^{-1}t +\Re(\gamma
t^{\lambda_2})+O(t^\kappa)$ and $\E \Xi = \gamma$ we obtain
 $\E W^{(t)}_r= O(t^\kappa)$. Since $W^{(t)}_r$ and $W^{(t)}_s$ are
 independent for $r\neq s$ conditionally on $(V_1,\ldots,V_b)$,
 it follows that
\begin{align}\label{eqd1}
\Delta(t)\le \Bigg\{\sum_{r=1}^b \E\Delta^2(V_rt)+
O\!\left(t^{2\kappa}\right) \Bigg\}^{1/2}+1,
\qquad t\ge1.
\end{align}
Now, we show that $\Delta(t)/t^\kappa=O(1)$. Note that this implies the
assertion. We denote
\begin{align*}
\Psi^\ast(t) := \sup_{1\le u\le t} \frac{\Delta(u)}{u^\kappa}.
\end{align*}
Then, (\ref{eqd1}) implies, that for appropriate $R>0$
\begin{align*}
\Psi^\ast(t)
\le
\Bigg\{\sum_{r=1}^b \E V_r^{2\kappa} (\Psi^\ast)^2(t)+R \Bigg\}^{1/2}+1,
\qquad t\ge1,
\end{align*}
and, with $\sqrt{a+b}\le \sqrt{a}+\sqrt{b}$ for $a,b\ge 0$ and
$\xi=\E \sum_{r=1}^b V_r^{2\kappa}<1$ this implies
\begin{align*}
\Psi^\ast(t)\le  \frac{\sqrt{R}+1}{1-\sqrt{\xi}}<\infty.
\end{align*}
The assertion follows.
%\end{proof}

\section{Examples} \label{Sex}

\begin{example}[Random splitting of intervals] \label{Ebinary}
Sibuya and Itoh \cite{SubItoh} studied the case of
random splitting of intervals, with uniformly distributed splitting
points; this is the case
$b=2$ and $\vv=(U,1-U)$, with $U\sim\unif(0,1)$.

We have
\begin{equation*}
  \phi(z)=\E U^z+\E(1-U)^z = 2\intoi u^z\dd u = \frac{2}{1+z},
\qquad \Re z>-1,
\end{equation*}
which is a rational function. The \che{}
\eqref{chareq} is $2/(1+\gl)=1$, and has the single
root $\gl=1$. Thus \refT{T1}(i) applies and shows
asymptotic normality, as stated by Dean and
Majumdar  \cite{DeanMaj}. Further,
$\ga=-\phi'(1)=1/2$, so \refT{TMV}(ii) yields
$\E\nx = m(x)=2x+O(x^\gd)$ for every $\gd>0$. More
precisely, \refT{Trational} yields
\begin{equation*}
  \E\nx=m(x)=2x-1, \qquad x\ge1,
\end{equation*}
which also can be shown directly from \eqref{a4} or from \eqref{sofie}.

For the asymptotic variance, we obtain from \refT{Trational}(ii),
since $\Nx=1$ and $a_0=-1$, using symmetry,
\begin{align*}
  \gb
&
=
\ga^{-1}\Bigpar{\E\bigpar{U+2U\bmin(1-U)+1-U}-2+1}
\\&\qquad{}
+2\ga^{-2}\Bigpar{2\E\bigpar{U(\ln(1-U)-\ln U)\ind{U<1-U}}}
-2\gai
\\&\qquad{}
+\ga^{-3}\Bigpar{2\E\bigpar{U\bmin(1-U)}}
-\gai
\\
&
= 20 \E\bigpar{U\bmin(1-U)}
+16\int_0^{1/2}u\bigpar{\ln(1-u)-\ln(u)}\dd u
-6
\\&
=8\ln 2-5
\approx
0.545177.
\end{align*}
This can also be obtained from \refT{TMV}(iii);
we have
\begin{align*}
  \psi(z,w)
&=
\E\Bigpar{(U^z+(1-U)^z)(U^w+(1-U)^w)}-\phi(z)\phi(w)
\\&
=
\frac2{1+z+w}+2B(z+1,w+1)-\frac4{(1+z)(1+w)}
\\&
=
\frac2{1+z+w}+2\frac{\Gamma(z+1)\Gamma(w+1)}{\Gamma(z+w+2)}-\frac4{(1+z)(1+w)}
\end{align*}
and thus
\begin{align*}
  \psi(1/2+\ii u,1/2-\ii u)
&=
1+\Gamma(3/2+\ii u)\Gamma(3/2-\ii u)-\frac4{|3/2+\ii u|^2}
\\&
=
1+|1/2+\ii u|^2\frac{\pi}{\cosh \pi u}-\frac4{|3/2+\ii u|^2},
\end{align*}
and, since $1-\phi(z)=(z-1)/(z+1)$,
\begin{align*}
\gb
=\frac1{\pi}\intooo
\Bigpar{
1+\frac{\pi}{\cosh \pi u}|1/2+\ii u|^2-\frac4{|3/2+\ii u|^2} }
\frac{|3/2+\ii u|^2}{|1/2+\ii u|^4}\dd u,
\end{align*}
which can be integrated (with some effort) to yield $8\ln 2-5$.

Consequently, by \refT{T1},
\begin{equation*}
 \frac{N(x)-2x}{\sqrt x}\dto \cN(0,8\ln2-5).
\end{equation*}
\end{example}

\begin{example}[$m$-ary splitting of intervals] \label{Emary}
We can generalize \refE{Ebinary} by splitting each
interval into $m$ parts, where $m\ge2$ is fixed,
using $m-1$ independent, uniformly distributed cut
points in each interval. This has been studied by
Dean and Majumdar  \cite{DeanMaj}.

We have $b=m$, and $V_1,\dots,V_m$ have the same distribution with
density $(m-1)(1-x)^{b-2}$, $0<x<1$.
Hence,
\begin{align*}
  \phi(z)
&= m\E V_1^z
=m(m-1)\intoi x^z(1-x)^z\dd x
= m(m-1)B(z+1,m-1)
\\&
= \frac{\Gamma(z+1)m!}{\Gamma(m+z)}
=\frac{m!}{(z+1)\dotsm(z+m-1)}.
\end{align*}
The characteristic equation $\phi(z)=1$ becomes
$\Gamma(z+m)/\Gamma(z+1)=m!$,
or
\begin{equation}
  \label{emary}
(z+1)\dotsm(z+m-1)=m!.
\end{equation}
The same equation appears in the analysis of
$m$-ary search trees.  It is shown by Mahmoud and
Pittel \cite{MahPit89}  and Fill and Kapur
\cite{FillKap03} that if $m\le 26$, then
$\Re\gl_2<1/2$, and thus (i) applies, but if
$m\ge27$, then $\Re\gl_2>1/2$, see also, e.g.,
Chauvin and Pouyanne \cite{chpo04} and Chern and
Hwang \cite{chhw01}. \refT{Trational} yields an
exact formula for $\E N(x)$ (although it is hardly
useful except when $m$ is small). It further leads
to a formula for the asymptotic variance, provided
$m\le26$.

We have, with $\psi(z)\=\Gamma'(z)/\Gamma(z)$ and
$H_z\=\psi(z+1)-\psi(1)$
(for integer $z$, these are the harmonic numbers)
\begin{equation*}
  \ga=-\phi'(1)=\psi(m+1)-\psi(2)=H_m-1.
\end{equation*}
\end{example}

\begin{example}[Random splitting of multidimensional intervals]
  \label{Equad}
Another generalization is to consider
$d$-dimensional intervals, where an interval is
split into $2^d$ subintervals by $d$ hyperplanes
orthogonal to the coordinate axis and passing
through a random, uniformly distributed point. This
too has been studied by Dean and Majumdar
\cite{DeanMaj}.

We have $b=2^d$. $V_1,\dots,V_b$ have the same distribution,
$V_j\eqd U_1\dotsm U_d$, where $U_k\sim U(0,1)$ are i.i.d.
Hence,
\begin{equation*}
  \phi(z)=2^d \E V_1^z = 2^d\bigpar{\E U_1^z}^d = \parfrac{2}{1+z}^d.
\end{equation*}
Again, $\phi$ is rational. The characteristic equation may be written
$\bigpar{(1+\gl)/2}^d=1$, with the roots
\begin{equation*}
  \gL=\set{2 e^{2\pi\ii k/d}-1: 0\le k\le d-1}.
\end{equation*}
Thus $\gs_2\=\Re\gl_2=2\cos\tfrac{2\pi}{d}-1$, and
the condition $\Re\gl_2<1/2$ is equivalent to
$\cos(2\pi/d)<3/4$, which holds for $d\le 8$, while
$\Re\gl_2>1/2$ for $d\ge9$. This justifies the
claims in Dean and Majumdar \cite{DeanMaj}.

The same characteristic equation, and the same phase transition,
appears for quad trees, see Chern, Fuchs and Hwang \cite{ChFHw}.

We further observe that $\ga=-\phi'(1)=d/2$.
\end{example}

The random trees in these three examples have also
been studied by \cite{ItohMah}, \cite{JavMahVah}
and \cite{JavVah}, where the properties of a
randomly selected branch are investigated. This
problem is quite different, and there is no phase
transition. See also \cite{SJ165}.

\begin{example}[Random splitting of simplices]
  \label{Esimplex}
Consider $d$-dimensional simplices,
where an simplex is split into $d+1$  new simplices by
choosing a random point $X$ in the interior and connecting it to the
vertices of the original simplex; each new simplex has as vertices $X$
and $d$ of the original $d+1$ vertices.

It is easily seen that this is equivalent to $d+1$-ary splitting as in
\refE{Emary}, see \cite[Lemma 3]{dev99}, so we have the same results
as there, with $m=d+1$.
In particular, $\nx$ is asymptotically normal if $d\le25$.
\end{example}

\begin{example}[Non-uniform splitting of intervals] \label{Ebeta}
Returning to binary splitting of intervals, we can generalize
\refE{Ebinary} by taking another distribution for the cut points;
we thus have
$b=2$ and $\vv=(V,1-V)$, where $V$ has any distribution on $(0,1)$.
An interesting case is when $V$ has a beta distribution
$V\sim B(a,a')$ with $a,a'>0$; then
\begin{equation*}
  \E V^z = B(a,a')\qi\intoi x^{z+a-1}(1-x)^{a'-1}\dd x
=\frac{B(a+z,a')}{B(a,a')}
=\frac{\Gamma(z+a)}{\Gamma(z+a+a')}\frac{\Gamma(a+a')}{\Gamma(a)};
\end{equation*}
$\E(1-V)^z$ is obtained by interchanging $a$ and $a'$.
In particular, if $a$ and $a'$ are integers, then $\phi$ is rational.

We consider two special cases.
\begin{romenumerate}
  \item
The symmetric case with $a'=a$, $V\sim B(a,a)$.
Then
\begin{equation*}
  \phi(z)=2\frac{\Gamma(z+a)}{\Gamma(z+2a)}\frac{\Gamma(2a)}{\Gamma(a)}
=\frac{\Gamma(z+a)}{\Gamma(z+2a)}\frac{\Gamma(1+2a)}{\Gamma(1+a)}.
\end{equation*}
We have
$%\begin{equation*}
  \ga=-\phi'(1)=H_{2a}-H_a,
$%\end{equation*}
with $H_x$ as in \refE{Emary}.
Numerical solution of the characteristic equation seems to show that
$\Re\gl_2<1/2$ if and only if $a<a_0$, where $a_0\approx 59.6$.

\item
The case $a'=1$, $V\sim B(a,1)$. Then
\begin{equation*}
\phi(z)
=\frac{\Gamma(z+a)}{\Gamma(z+a+1)}\frac{\Gamma(a+1)}{\Gamma(a)}
+\frac{\Gamma(z+1)}{\Gamma(z+a+1)}\Gamma(a+1)
=\frac{a}{z+a}+\frac{\Gamma(z+1)}{\Gamma(z+a+1)}\Gamma(a+1).
\end{equation*}
One finds   $\ga=H_{a}/(a+1)$.
The characteristic equation $\phi(\gl)=1$ is equivalent to
$\Gamma(a+1)\Gamma(\gl+1)/\Gamma(\gl+a+1)=\xfrac{\gl}{(\gl+a)}$
or
\begin{equation*}
  \frac{\Gamma(a+\gl)}{\Gamma(\gl)}=\Gamma(a+1).
\end{equation*}
When $a=m$ is an integer, this is the same as
\eqref{emary}, so  $\Re\gl_2<1/2$ for integer $a$
if and only if $a\le26$. In general, numerical
solution of the characteristic equation seems to
show that $\Re\gl_2<1/2$ if and only if $a<a_0$,
where $a_0\approx 26.9$.
\end{romenumerate}
\end{example}

\section{Non-examples}\label{Snex}
In this section, we give a few examples where our theorems are not valid.

\begin{example}[Lattice]
  \label{Elattice}
In the lattice case, there exists $R>1$ such that every
$V_j\in\set{R^{-k}:k\ge1}\cup\set0$ a.s.
In this case, $\phi$ is periodic with period $2\pi\ii/\ln R$; in
particular, the \che{} \eqref{chareq} has infinitely many roots
$1+2\pi\ii n/\ln R$ on \set{\gl:\Re\gl=1}, and thus \conditionbd{}
fails.
Indeed, it is obvious from \eqref{a4} that $\nx=N(R^m)$ when $R^m\le
x<R^{m+1}$, so $\E\nx/x$ oscillates and does not converge as \xtoo.
The natural approach is to consider only $x\in\set{R^m:m\ge0}$.
It is then straightforward to prove an analogue of \refT{TMV}, using
the lattice versions of the renewal theory theorems that were used in
\refS{Smean}. An analogue of \refT{T1} then follows by the usual
(discrete) contraction method, as in \cite{neru4}.
We leave the details to the reader.
\end{example}

\begin{example}[Deterministic]
  \label{Edeterm}
If $V=(V_1,\dots,V_b)$ is deterministic, then so is $N(x)$, and it is
meaningless to ask for an asymptotic distribution. However, it makes
sense to study the asymptotics of $N(x)=m(x)$. (Clearly, $\gss(x)=0$.)

If $V$ is non-lattice, then $N(x)/x\to\ga$ by
\refT{TMV} and \refR{Rnl}.  If $V$ is lattice, we
consider, as in \refE{Elattice}, only $x=R^m$,
$m\ge1$.

We may assume that $V_j>0$ for each $j$.
By
Dirichlet's theorem \cite[Theorem 201]{HW},
for every $\eps>0$, there exist arbitrarily large $t$ such that
$| V_j^{\ii t}-1| <\eps$ for $j=1,\dots,b$;
thus $\limsup_{\ttoo}|\phi(1+\ii t)|=1$.
Hence \condition{B(1)} does not hold, and therefore, by \refL{LBd},
\conditionbd{} does not hold for any $\gd\le1$.

More precisely, if $| V_j^{\ii t}-1| <\eps$ for $j=1,\dots,b$, let
$z_0=1+\ii t$.
Then $|\phi(z_0)-1|<\eps$ and
$$
|\phi'(z_0)+\ga|
=\left|\sumjb\ln V_j (V_j^{1+it}-V_j)\right|
\le\eps\ga.
$$
Since further $|\phi''(z)|\le\sum_j|\ln V_j|^2$ for $\Re z\ge0$,
it follows easily that if $\eps$ is small enough, then $\phi(z)-1$ has
a zero in the disc $B\=\set{z:|z-z_0|<2\eps/\ga}$.
(Use the Newton--Raphson method, or Rouch\'e's theorem and a comparison
with the linear function $\phi(z_0)+(z-z_0)\phi'(z_0)$.)
It follows that there exists a sequence $\gl_n\in\gL$ with
$\Re\gl_n\to1$ and $\Im\gl_n\to+\infty$.

We give some concrete examples:

$V=(1/2,1/2)$ is lattice with $R=2$ and $N(2^n)=2^n$.

$V=(\tau\qi,\tau^{-2})$ where $\tau=(1+\sqrt5)/2$ (the golden ratio) is
lattice with $R=\tau$ and $N(\tau^n)=F_{n+3}-1$, $n\ge0$, as is easily
proven by induction. ($F_n$ denotes the Fibonacci numbers.)
Thus, $N(\tau^n)\sim 5\qhi\tau^{n+3}$.

$V=(1/3,2/3)$ is non-lattice and thus $N(x)\sim\ga\qi x$, where
$\ga=\tfrac13\ln 3+\tfrac23\ln(3/2)=\ln3-\tfrac23\ln2$.
\end{example}

\section{Some related models}\label{Srelated}

The basic model may be varied in various ways. We mention here some
variations that we find interesting.
We do \emph{not} consider these versions in the present paper;
we leave the
possibility of extensions
of our results as an open problem, hoping that these remarks will
be an inspiration for future research.

\begin{remark}
  By our assumptions, the label of a node equals the sum of the labels
  of its children. Another version would be to allow a (possibly
  random) loss at each node.
One important case is \emph{\Renyi's parking problem}
\cite{Renyi}, where a node with label $x$ is interpreted as an
interval of length $x$ on a street, where cars of length 1 park at
random.
Each car splits an interval of length $x\ge1$ into two free intervals with
  the lengths $U(x-1)$ and $(1-U)(x-1)$, where $U\sim \unif(0,1)$.
An obvious generalization is to split $(x-1)$ using an arbitrary  random vector
\vvb.
(The one-sided version, where we study only one branch of the tree, is
  studied in \cite{ItohMah}, \cite{SJ165}.)
\end{remark}

\begin{remark}\label{Rrandfrag}
Krapivsky, Ben-Naim and Grosse \cite{KBG01,KBG04} have studied
a fragmentation process where fragmentation stops stochastically, with
a probability $p(x)$ of further fragmentation that in general depends
on the size $x$ of the fragment. Our process is the case $p(x)=\ind{x\ge1}$.
Another interesting case is $p(x)=1-e^{-x}$, see \refR{Rdiscrete} below.
\end{remark}

\begin{remark}\label{Rdiscrete}
Our model is a continuous version of the split trees studied by
  Devroye \cite{dev99}, where the labels are integers (interpreted as
numbers of balls to be distributed in the corresponding subtree)
and each label $n$ is, except at the leaves, randomly split according
  to a certain procedure
into $b$ integers summing to $n-s_0$; here $s_0$ is a small positive
  integer
(for example 1) that represents the number of balls stored at the node.
Typical examples are binary search trees, $m$-ary search trees and
  quadtrees.
We can regard the continuous model as an
approximation of the discrete, or conversely, and
it is easy to guess that many properties will have
similar asymptotics for the two models. This has
been observed in several examples by various
authors, see \cite{DeanMaj} and \cite{ChFHw}. For
example, the results for \refE{Emary} parallel
those found for $m$-ary search trees by
\cite{MahPit89}, \cite{chhw01}, \cite{FillKap03},
\cite{chpo04} and others. Similarly, the results in
\refE{Equad} parallel those found for
  quadtrees by \cite{ChFHw}.

We study only the continuous version in this paper.
It would be very interesting to be able to rigorously transfer results
from the continuous to the discrete version (or conversely); we will,
however, not attempt this here.

Note that for binary search trees, we have $n$ random (uniformly
distributed) points in an
interval, split the interval by the first of these points, and
continue recursively splitting each subinterval that contains at least
one of the points. If we scale the initial interval to have length
$n$, then the probability that a subinterval of length $x$ contains at
least one point is $\approx1-e^{-x}$. Thus it seems likely that the
binary search tree is well approximated by a fragmentation tree, with
\vv{} as in \refE{Ebinary}, with a fragmentation probability
$1-e^{-x}$ as in \refR{Rrandfrag}. The same goes for random quadtrees and
simplex trees corresponding to Examples \refand{Equad}{Esimplex}.
\end{remark}

\newcommand\jour{\emph}
\newcommand\book{\emph}
\newcommand\inbook{\emph}
\newcommand\vol{\emph}
\newcommand\no[1]{\unskip}

\newcommand\toappear{\unskip, to appear}
\newcommand\webcite[1]{\hfil\penalty0\texttt{\def~{\~{}}#1}\hfill\hfill}
\newcommand\webcitesvante{\webcite{http://www.math.uu.se/\~{}svante/papers}}

\end{document}